\documentclass[AMA,STIX1COL]{WileyNJD-v2}
\usepackage{moreverb}

\newcommand\BibTeX{{\rmfamily B\kern-.05em \textsc{i\kern-.025em b}\kern-.08em
T\kern-.1667em\lower.7ex\hbox{E}\kern-.125emX}}

\articletype{Original Article}

\received{<day> <Month>, <year>}
\revised{<day> <Month>, <year>}
\accepted{<day> <Month>, <year>}

\usepackage{subfiles}
\usepackage{nicefrac}
\usepackage{amsmath,amsthm,amsfonts,mathtools}
\usepackage{graphicx}
\usepackage{float}
\usepackage{textcomp}
\usepackage{caption}
\usepackage{multirow}
\usepackage{subcaption}
\usepackage{hyperref}
\hypersetup{
    colorlinks=true,
    linkcolor=blue,
    filecolor=magenta,      
    urlcolor=cyan,
}
\usepackage{xcolor}

\begin{document}


\title{Reference Governor for Constrained Spacecraft Orbital Transfers}

\author[1]{Simone Semeraro}

\author[1]{Ilya Kolmanovsky}

\author[2]{Emanuele Garone}

\authormark{SEMERARO \textsc{et al.}}

\address[1]{\orgdiv{Department of Aerospace Engineering}, \orgname{The University of Michigan}, \orgaddress{\city{Ann Arbor}, \state{Michigan}, \country{USA}}}

\address[2]{\orgdiv{École Polytechnique de Bruxelles}, \orgname{Université Libre de Bruxelles (ULB)}, \orgaddress{\city{Brussels}, \country{Belgium}}}


\corres{Ilya Kolmanovsky, Department of Aerospace Engineering, The University of Michigan, 1320 Beal Ave, Ann Arbor, MI 48109, USA. \email{ilya@umich.edu}}

\fundingInfo{National Science Foundation Grant Numbers CMMI-1904394 and ECCS-1931738.}


\abstract[Abstract]{
The paper considers the application of feedback control to orbital transfer maneuvers subject to constraints on the spacecraft thrust and on avoiding the collision with the primary body.
Incremental reference governor (IRG) strategies are developed to complement the nominal Lyapunov controller, derived based on Gauss Variational Equations, and enforce the constraints.  Simulation results are reported that demonstrate the successful constrained orbital transfer maneuvers with the proposed approach. A Lyapunov function based IRG and a prediction-based IRG are compared. While both implementations successfully enforce the constraints, a prediction-based IRG is shown to result in faster maneuvers.

}

\keywords{Gauss Variational Equations, Constrained Control, Nonlinear Control, Reference Governors}

\jnlcitation{\cname{
\author{S. Semeraro}, 
\author{I. Kolmanovsky}, and 
\author{E. Garone}} (\cyear{2022}), 
\ctitle{Reference Governor for Constrained Spacecraft Orbital Transfers}, \cjournal{Advanced Control for Applications}, \cvol{}.}

\maketitle

\section{Introduction}\label{sec:intro}

In this paper we consider the problem of spacecraft orbital transfer in a setting of a Two Body problem \cite{battin1999introduction,gurfil2016celestial}. Our approach is distinguished by using the reference governor\cite{garone2017reference} to complement the nominal Lyapunov feedback controller and enforce pointwise-in-time state and control constraints.  Such constraints in orbital transfer problems could be imposed to   maintain sufficient separation distance to the planet (e.g., to avoid planetary collisions), to avoid the exceedance of the maximum thrust limit, and to maintain trajectories in model validity range.  

Orbital transfers are critical components of many space missions. For instance, an Earth orbiting spacecraft is often launched into a parking orbit first and then must execute a transfer maneuver into a Geostationary Equatorial Orbit (GEO), which is the primary destination for many communication and Earth observation satellites. A spacecraft flying to Mars has to execute an Earth to Mars orbital transfer around the Sun.

Traditionally, the design of orbital transfer trajectories has been approached using the theory of optimal control, where a typical objective is minimizing the cumulative fuel consumption. Once an optimal open-loop trajectory has been determined, several Trajectory Correction Maneuver (TCM) points are defined at which orbital corrections could be made.  

The feedback control for orbital transfers has also been investigated \cite{gurfil2016celestial,petropoulos2005refinements,holt2020low,chang2002lyapunov}, in particular, for missions that exploit the low thrust propulsion.  
For instance, the Q-law \cite{petropoulos2004low,petropoulos2005refinements} has been  proposed and extensively studied. 

The advantages of feedback control laws are that they apply to continuously thrusting systems and provide a way to apply systematic corrections all the time as opposed to only at pre-specified TCM points.  As in other applications, the use of feedback for orbital transfers is appealing due to improved robustness. 

At the same time, nonlinear dynamics of spacecraft motion and constraints, such as limits on the thrust forces that are much smaller as compared to the gravitational forces, complicate the development of effective feedback solutions.  Many techniques, such as feedback linearization, while formally applicable to spacecraft equations of motion, may not be practical due to large thrust and fuel consumption they require.

In this paper we consider the application of reference governors to assist Lyapunov-based feedback controllers in satisfying pointwise-in-time state and control constraints during orbital transfer maneuvers. Reference governors \cite{garone2017reference} are add-on schemes to nominal controllers that minimally modify reference commands (set-points) when it becomes necessary to satisfy the constraints.  For the nominal controller in our orbital transfer problem, we choose a Lyapunov-based controller\cite{gurfil2016celestial}, derived based on Gauss Variational Equations (GVEs), for which we highlight the connection between closed-loop stability and persistence of excitation conditions in achieving convergence to the target orbit; this connection is not made explicit in the prior literature \cite{gurfil2016celestial,hattenthesis}.  

We then propose a modification of an incremental reference governor
(IRG)  strategy\cite{tsourapas2009incremental}  that we refer to as a multi-step multi-mode incremental reference governor (m$^2$IRG) to handle constraints on minimum radius of periapsis, thrust and eccentricity.  The m$^2$IRG checks the safety of changing the reference (in our case, the reference will be defined by $5$ orbital elements of the target orbit) and either increments the reference towards the ultimate target (if safe) or maintains the current reference (if unsafe).
The m$^2$IRG is different from IRG\cite{tsourapas2009incremental} in that it takes advantage of time available for computations during the orbital transfer to check safety of multiple increments of the reference. In this regard, it is closer to nonlinear command governor in which optimization is implemented inexactly\cite{garone2022command}.
Additionally,  m$^2$IRG can exploit switching between a set of feedback gains (corresponding to multiple modes) for the nominal Lyapunov controller.  

Two implementations of m$^2$IRG are considered and compared that are different in the way safety of an incremented reference is checked. The first implementation exploits the knowledge of the Lyapunov function for the closed-loop system and invariance of its sublevel sets; hence it checks whether the current state is in a safe sublevel set of the Lyapunov function to check the safety of the incremented reference.  This check requires that maximum of each constraint over a sublevel set around the incremented reference be determined, and reduces to solving several low dimensional nonlinear optimization problems.  For this Lyapunov function based  m$^2$IRG   switching between multiple feedback gains can help better align the sublevel sets of the closed-loop Lyapunov function with the geometry of the constraints which enables larger reference increments.  The second implementation checks constraints based on the online prediction of the closed-loop response to incremented reference over a long prediction horizon.  

The paper is organized as follows: Section~\ref{sec:2} presents the equations of motion in the form of GVEs. In Section~\ref{sec:Lyapunov cont} we describe the nominal Lyapunov controller and highlight the need for persistence of excitation conditions to ensure closed-loop stability. Section~\ref{sec:IRG} introduces constraints as well as m$^2$IRG  to handle them. Section~\ref{sec:5} reports the results of applying m$^2$IRG to simulated orbital transfers. 
Section~\ref{sec:7} demonstrates the robustness of the proposed solution to changes in the mass of the spacecraft during the maneuver that leads to constraints becoming time-varying.  Section~\ref{sec:8} introduces an online prediction-based  m$^2$IRG and reports the corresponding simulation results. It is shown that the online prediction-based m$^2$IRG leads to faster maneuvers as compared to the Lyapunov function based m$^2$IRG. Finally, Section~\ref{sec:9} provides a summary and concluding remarks.

\section{Gauss-Euler Variational Equations}\label{sec:2}
The spacecraft translational dynamics in orbital transfer problems can be represented 
by different sets of equations.  A common approach relies on the direct application 
of 2nd Newton's law leading to six differential equations of motion (EoMs) for three components of spacecraft position vector, $\vec{r}$, and three components of spacecraft 
velocity vector, $\vec{v}$ based on
\begin{equation}\label{equ:newton}
\begin{array}{c}
    \displaystyle \dot{\vec{r}}=\vec{v} \\
    \displaystyle \dot{\vec{v}}=-\frac{\mu}{r^3} \vec{r}+\frac{\vec{F}}{m},
\end{array}
\end{equation}
where $\mu$ is the gravitational parameter of the primary body
($398600.436$  km$^3$s$^{-2}$ for Earth), $\vec{F}$ is the net  ``perturbation'' force (i.e., the net of all forces except for the ideal inverse distance squared gravity) applied to the spacecraft, and $m$ is the spacecraft mass.

\begin{figure}[h]
\centering
\includegraphics[width=8cm]{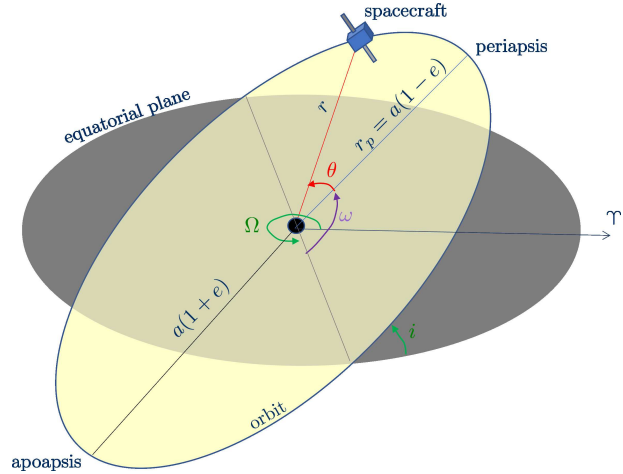}
\caption{Classical orbital elements.}\label{fig:basicoe}
\end{figure}

An alternative approach is the use of Gauss-Euler Variational Equations (GVEs) \cite{gurfil2016celestial,battin1999introduction} derived by the method of variation of constants of Lagrange.
Such a formulation of EoMs applies to the six classical orbital elements of Keplerian two body 
problem in presence of perturbation forces such as spacecraft thrust.
These orbital elements are  the semi-major axis $a$ [km], the eccentricity $e$,  the inclination $i$ [rad],  the Right Ascension of the Ascending node $\Omega$ [rad], the argument of periapsis $\omega$ [rad] and the spacecraft true anomaly $\theta$ [rad].
See Figure~\ref{fig:basicoe}.
The orbital elements $a$ and $e$ determine the size and shape of the orbit, the orbital elements $i$, $\Omega$ and $\omega$ determine the orbit orientation in an inertial frame
and $\theta$ determines the current spacecraft position within the orbit.  Given the values of the six orbital elements one can uniquely reconstruct $\vec{r}$ and $\vec{v}$.  The inverse transformation is not defined everywhere, and this limits the applicability of GVEs to  orbits with $a>0$, $e>0$ and
$i \neq 0$.  Variants of GVEs can be derived for nonclassical equinoctual orbital elements \cite{battin1999introduction}
that can be used for zero inclination or eccentricity orbits; we leave the treatment of such cases to future publications.

To model the evolution of the orbital elements, a moving STW frame is introduced with the origin at the center of mass of spacecraft and with the unit vectors $\hat{e}_r$, $\hat{e}_\theta$ and $\hat{e}_h$
defined according to
$$\hat{e}_r=\frac{\vec{r}}{r},~\hat{e}_h=\frac{\vec{h}}{h},~ \vec{h} =\vec{r} \times \vec{v},~
\mbox{and}~ \hat{e}_\theta=\hat{e}_h \times \hat{e}_r,$$
where $\times$ denotes the vector product.
Decomposing the thrust force per unit mass 
[km/s$^2$]  applied to the spacecraft as
\begin{align*}
    \frac{\vec{F}}{m} = S\hat{e}_r + T\hat{e}_\theta + W\hat{e}_h,
\end{align*}
where $m$ denotes the mass of the spacecraft,
the GVEs take the following form:
\begin{align}\label{gve}
    \frac{da}{dt} &= \frac{2a^2}{\sqrt{\mu p}}e\sin{\theta}S+\frac{2a^2}{\sqrt{\mu p}}\frac{p}{r}T, \nonumber \\
    \frac{de}{dt} &= \frac{p\sin\theta}{\sqrt{\mu p}}S + \frac{p(\cos\psi+\cos\theta)}{\sqrt{\mu p}}T, \nonumber\\
    \frac{di}{dt} &= \frac{r}{\sqrt{\mu p}}\cos{(\theta+\omega)}W, \\
    \frac{d\Omega}{dt} &= \frac{r}{\sqrt{\mu p}}\frac{\sin{(\theta+\omega)}}{\sin i}W, \nonumber\\
    \frac{d\omega}{dt} &= -\frac{p\cos\theta}{e\sqrt{\mu p}}S+\frac{(r+p)\sin\theta}{e\sqrt{\mu p}}T-\frac{r\sin(\theta+\omega)\cot i}{\sqrt{\mu p}}W,\nonumber \\
    \frac{d\theta}{dt} &= \frac{\sqrt{\mu p}}{r^2}+\frac{p\cos\theta}{e}\frac{S}{\sqrt{\mu p}}-\frac{p+r}{e}\cos\theta\frac{T}{\sqrt{\mu p}}, \nonumber
\end{align}
where  $r=\frac{p}{1+e\cos\theta}$ is the distance from the gravity center to the spacecraft center of mass, $p=a(1-e^2)$ is the orbit parameter (also called semi-latus rectum), and $\psi = \arccos\left(\frac{1}{e}-\frac{r}{ae}\right)$ is the eccentric anomaly.

Note that with $S=T=W=0$, the trajectories of $a$, $e$, $i$, $\Omega$, $\omega$ are constants in time while $\theta$ will evolve as a function of time corresponding to the motion of the spacecraft on an elliptic orbit under the ideal force of gravity only.  Note also that the unit vectors $\hat{e}_r$, $\hat{e}_\theta$ and $\hat{e}_h$ can be computed based on the current position and velocity vectors of the spacecraft; thus a thrust vector specified on the basis of a feedback law for $S$, $T$ and $W$ components
(i.e., in STW frame) can be transformed into an inertial frame using the computed onboard Direction Cosine Matrix (DCM) between the two frames.

\section{Lyapunov Controller and Need for Persistence of Excitation}\label{sec:Lyapunov cont}
The reference governor augments a nominal controller which in this paper is a Lyapunov controller of Jurdjevic-Quinn damping feedback type \cite{sepulchre2012constructive}.
Let
\begin{align*}
    X = [a \ e \ i \ \Omega \ \omega]^{\sf T}, \ \ U = [S \ T \ W]^{\sf T}.
\end{align*}
Then the first five of GVEs (\ref{gve}) can be written in the following condensed form:
\begin{equation} \label{equ:main1}
    \dot X(t) = G\big(X(t),\theta(t)\big)U(t),
\end{equation}
where $G(X(t),\theta(t)) \in \mathbb{R}^{5 \times 3}$.
Note that since the objective is to steer the spacecraft into a target orbit but not to a particular location in that orbit (each location, once in the target orbit, will be visited eventually due to periodicity of the orbit) the true anomaly $\theta(t)$ is not included into the state vector $X(t)$ and is treated as a time-varying parameter, the evolution of which is determined by the sixth equation
in (\ref{gve}).  With a slight abuse of notation, and in order to shorten the expressions, in the sequel we designate
$$
G(t) = G(X(t),\theta(t)).
$$

The dynamics in (\ref{equ:main1}) are drift-free.  Hence any point in $\mathbb{R}^5$ in the range of definition of (\ref{gve}) is an unforced equilibrium (i.e., an equilibrium with $U=0$) and can be commanded as the target.  Despite the system having only three inputs and five states, the dynamics are small time locally controllable. This can be concluded from the original EoMs (\ref{equ:newton}) that are based on Newton's law and can be feedback linearized into a system of three uncoupled double integrators, hence clearly controllable. 

There is often a topological obstruction to stabilizability of drift-free systems with more states than inputs due to Brockett's necessary condition \cite{brockett1983asymptotic}.  However, 
(\ref{equ:main1}) is not a time-invariant system and turns out to be stabilizable. 

Let $P \succ 0$ be a $5 \times 5$ positive-definite matrix and let
\begin{equation}\label{equ:main2}
    \Tilde{E}(t)=P\left(X(t)-X_{\tt des}\right),
\end{equation}
denote the weighted tracking error where $X_{\tt des}$ is the vector of five orbital elements of the target orbit.

We choose a Lyapunov function candidate as 
\begin{equation}\label{equ:main2.5}
V(\Tilde{E})=\frac{1}{2}\Tilde{E}^{\sf T}(t)P^{-1}\Tilde{E}(t).
\end{equation} 
Its time derivative along the trajectories of (\ref{equ:main1}) is given by
\begin{equation*}
    \dot V = \frac{\partial V}{\partial\Tilde{E}}\dot{\Tilde{E}}=\Tilde{E}^{\sf T}P^{-1}\dot{\Tilde{E}}=\Tilde{E}^{\sf T}P^{-1} P G(t) U=\Tilde{E}^{\sf T} G(t) U.
\end{equation*}
Defining a feedback law as
\begin{equation}\label{equ:main3}
U = -G^{\sf T}(t)  \Tilde{E}=  -G^{\sf T}(t)P\left(X(t)-X_{\tt des}\right) ,
\end{equation}
we obtain,
\begin{equation}\label{equ:main4.5}
    \dot V = -\Tilde{E}^{\sf T} G(t) G^{\sf T}(t) \Tilde{E}.
\end{equation}
While $\dot{V} \leq 0$, implying the invariance of sublevel sets of $V$, it is not sufficient to prove the convergence of $\Tilde{E}(t)$ to zero as $t \to \infty$.  Using the LaSalle's invariance principle as attempted in prior work \cite{gurfil2016celestial} is complicated by the fact that the system is time-varying and that the matrix $G(t) G^{\sf T}(t)$ is $5 \times 5$ but is at most of rank $3$ as $G(t)$ is $3 \times 5$. Hence it is not immediately clear how to show that there are no closed-loop trajectories, except for $\tilde{E}=0$, that could evolve in the time-varying null space of $G(t) G^{\sf T}(t)$. 

The desired result can be obtained under a ``persistence of excitation'' assumption following the same lines of proof as in adaptive control\cite{tao2003adapt}. For completeness, given a slightly different form with the $P$ matrix in (\ref{equ:main4}), and since the decay estimates coul be exploited in the construction and anaysis of the reference governor we show the details. First, note that the closed-loop dynamics are linear time-varying and given by
\begin{equation}\label{equ:main4}
\dot{\Tilde{E}}=-PG(t)G^{\sf T}(t) \tilde{E},
\end{equation}
and we let $\Phi(t_0,t)$ denote the state transition matrix
of the linear time-varying system (\ref{equ:main4}). Consider the change in the Lyapunov function between time instants $\sigma_0$ and $\sigma_0+\delta_0>\sigma_0$.  We have
\begin{align*}
    &V(\Tilde{E}(\sigma_0+\delta_0)) = V(\Tilde{E}(\sigma_0)) + \int_{\sigma_0}^{\sigma_0+\delta_0} \dot V(\tilde{E}(\tau))d\tau 
    = V(\Tilde{E}(\sigma_0))-\int_{\sigma_0}^{\sigma_0+\delta_0} \Tilde{E}^{\sf T}(\tau)G(\tau)G^{\sf T}(\tau)\Tilde{E}(\tau)d\tau = \\
    = &V(\Tilde{E}(\sigma_0)) - \Tilde{E}^{\sf T}(\sigma_0) \int_{\sigma_0}^{\sigma_0+\delta_0}\Phi^{\sf T}(\tau,\sigma_0)G(\tau)G^{\sf T}(\tau)\Phi(\tau,\sigma_0) d\tau\Tilde{E}(\sigma_0) 
    = V(\Tilde{E}(\sigma_0)) - \Tilde{E}^{\sf T}(\sigma_0) \tilde{O}(\sigma_0,\sigma_0+\delta_0)  \Tilde{E}(\sigma_0),
\end{align*}
where $\tilde{O}(\sigma_0,\sigma_0+\delta_0)$ designates the observability grammian for the system, 
\begin{align}\label{equ:main5}
& \dot{\Tilde{E}}=-PG(t)G^{\sf T}(t) \tilde{E}, \\
& y(t)=G^{\sf T}(t) \Tilde{E}. \nonumber
\end{align}
If there exists $\tilde{\epsilon}>0$ and $\delta_0>0$ such that for all $\sigma_0>0$,
\begin{equation}\label{equ:main6}
\tilde{O}(\sigma_0,\sigma_0+\delta_0) \geq \tilde{\epsilon} I_{5 \times 5},
\end{equation}
where $I_{5 \times 5}$ is the $5\times 5$ identity matrix; then
$$ V(\Tilde{E}(\sigma_0+\delta_0)) \leq V(\Tilde{E}(\sigma_0))
- \tilde{\epsilon} \Tilde{E}^{\sf T}(\sigma_0)   \Tilde{E}(\sigma_0)  \leq  V(\Tilde{E}(\sigma_0))
- \frac{2\tilde{\epsilon}}{\lambda_{\tt max}(P^{-1})} V(\Tilde{E}(\sigma_0))
=\left(1- 2\tilde{\epsilon}\lambda_{\tt min}(P)\right) V(\Tilde{E}(\sigma_0))
,$$
where we used (\ref{equ:main2}) in the second to last step.
Applying this recursively leads to 
\begin{equation}\label{equ:L1}
V(\Tilde{E}(\sigma_0+k \delta_0)) \leq
\left( 1- 2 \tilde{\epsilon}\lambda_{\tt min}(P)\right)^k V(\Tilde{E}(\sigma_0)).
\end{equation}
Note that $ 1- 2 \tilde{\epsilon}\lambda_{\tt min}(P)<1$, while
$ 1- 2 \tilde{\epsilon}\lambda_{\tt min}(P)<0$ and $V(\Tilde{E}(\sigma_0))>0$ would imply $V(\Tilde{E}(\sigma_0+ \delta_0))<0$, which is not possible since $P^{-1} \succ 0$
in (\ref{equ:main2}).
Hence $V(\Tilde{E}(\sigma_0+k \delta_0)) \to 0$ and  $\Tilde{E}(t_0+k\delta_0) \to 0$ 
as $k \to \infty$, where $k$ is an integer.  For $t \in (t_0+k \delta_0,t_0+(k+1)\delta_0]$,
\begin{equation}\label{equ:L2}
\frac{1}{2}\lambda_{\tt min}(P^{-1}) \|\tilde{E}(t)\|^2_2 \leq V(\tilde{E}(t)) \leq V(\tilde{E}(t_0+k \delta_0)), 
\end{equation}
and hence $V(\tilde{E}(t)) \to 0$ and $\tilde{E}(t) \to 0$ as $t \to \infty$. The remaining step is to simplify (\ref{equ:main6}). As the property of uniform complete observability is invariant under output feedback\cite{tao2003adapt},  a similar property to
(\ref{equ:main6}) must hold for 
\begin{align}\label{equ:main7}
& \dot{\Tilde{E}}=
-PG(t)G^{\sf T}(t) \tilde{E}
+PG(t)y(t)=0, \\
& y(t)=G^{\sf T}(t) \Tilde{E}, \nonumber
\end{align}
i.e., we can ensure that (\ref{equ:main6}) holds by showing that there exists ${\epsilon}>0$ and $\delta_0>0$ such that
for all $\sigma_0>0$,
\begin{equation}\label{equ:main6mod}
{O}(\sigma_0,\sigma_0+\delta_0)
= \int_{\sigma_0}^{\sigma_0+\delta_0} G(\tau) G^{\sf T}(\tau) d \tau
\geq {\epsilon} I_{5 \times 5}.
\end{equation}
The condition (\ref{equ:main6mod}) is the persistence of excitation condition on the ``regressor'' matrix, $G(t)$. It is possible\cite{zhang2015gramian} to obtain a relation between ${\epsilon}$ in (\ref{equ:main6mod}) and $\tilde{\epsilon}$ in (\ref{equ:main6}) which determines the convergence rate of $V(\tilde{E}(t))$ to zero.

\begin{figure}[h]
    \centering
    \includegraphics[width = 6.5 cm]{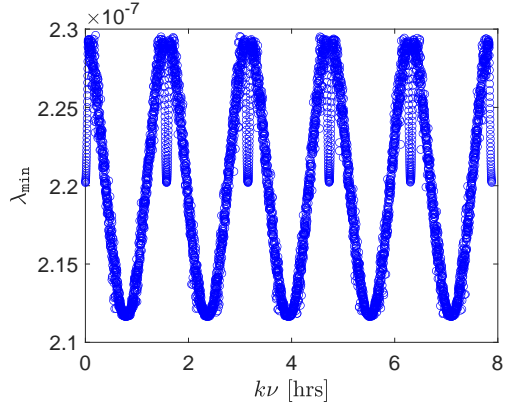}
    \caption{Minimum eigenvalue of 
    ${O}(k\nu, (k+1)\nu))$, $\nu=5$ sec versus time. }
    \label{fig:timeintPE}
\end{figure}

Verifying (\ref{equ:main6mod}) analytically for (\ref{equ:main1}), (\ref{equ:main3}) has been attempted through symbolic computations but does not appear to be tractable even for the case of the reduced order subsystem with states $a,e,\omega,\theta$ and control inputs $S$ and $T$ (i.e., $W=0$). However, (\ref{equ:main6mod}) can be checked by numerical integration. One approach is to choose $\nu>0$ small and compute 
$ \lambda_{\tt min}({O}(k\nu, (k+1)\nu ))$ for $k=0,1,2,\cdots$
over a given orbit or trajectory.   Noting that, by Weyl's inequality, 
$$\lambda_{min}({O}(k\nu, (k+m)\nu)) \geq \sum_{j=0}^{m-1}
\lambda_{min}({O}(k\nu + j\nu, k\nu+ (j+1)\nu  )),$$
the values of $\delta_0$ and $\epsilon$ that verify (\ref{equ:main6mod}) can be identified.  Using this procedure for specific orbits and transient maneuvers in this paper, we have not found even a single case where (\ref{equ:main6mod}) does not hold. In Figure~\ref{fig:timeintPE}, we consider a (periodic) orbit with $a=6878$, $e=0.02$, $i=\frac{\pi}{2},$ 
$\Omega=\frac{3\pi}{2},$
$\omega=\pi$ and plot 
$\lambda_{\tt min}(\bar{O}(k\nu, (k+1)\nu )),$ where $\nu=5$ sec.  Then, for instance, $\delta_0 = 180$ s, $\epsilon$ = 2$\cdot10^{-7}$ verify condition (\ref{equ:main6mod}).  By continuity, we also expect (\ref{equ:main6mod}) to hold for trajectories in a neighborhood of this nominal orbit.

Note that (\ref{equ:main3}) can be viewed as an approximate solution to an unconstrained minimization problem of  a regularized variant of $V(\tilde{E}(t+\Delta t))$, for small $\Delta t>0$ with respect to $U(t)$, i.e., of
$$V(\tilde{E}(t+\Delta t))+\frac{1}{2} \int_t^{t+\Delta t} U^{\sf T}(\tau) U(\tau) d \tau \approx V(\tilde{E}(t))+\dot{V}(t) \Delta t+\frac{1}{2} U^{\sf T}(t) U(t) \Delta t
=V(\tilde{E}(t))+\tilde{E}^{\sf T}(t) G(t) U(t) \Delta t+\frac{1}{2} U^{\sf T}(t) U(t) \Delta t.
$$
Since $V(\tilde{E})=\frac{1}{2} \tilde{E}^{\sf T} P^{-1} \tilde{E}=\frac{1}{2} (X-X_{\tt des})^{\sf T} P (X-X_{\tt des})$, the matrix $P$ in the controller calibration phase can be tuned to assign higher penalty to those orbital elements the response of which needs to be sped up in a similar manner as state weighting matrices are employed in Linear Quadratic or Model Predictive Control.

\section{Constraints and Lyapunov Function Based Incremental Reference Governor}\label{sec:IRG}
\subsection{Constraints}

Three constraints are imposed on the spacecraft translational motion during the orbital transfer maneuver.  

The first constraint has the form,
\begin{equation}\label{equ:maincon1}
    c_1(X) = r_p - r_{\tt min} \geq 0,\quad r_p=a(1-e). 
\end{equation}
It ensures that the radius of the periapsis of the spacecraft orbit, $r_p$
is larger than $r_{\tt min}$.  The orbital elements which evolve
according to GVEs (\ref{gve}) are that of the osculating orbit, that is of the Keplerian orbit which results if the perturbation force components are set to zero, i.e., $S=T=W=0$.  Thus the constraint (\ref{equ:maincon1})
protects not only against the distance to the primary $r$ falling below $r_{\tt min}$ at any given time instant and thus avoiding being too close to/colliding with the primary but also that $r$ will stay above $r_{\tt min}$ even if there is a thruster failure and thrust becomes zero.

The second constraint is imposed on the spacecraft thrust magnitude not to exceed the values that the thruster can actually deliver.  This constraint has the following form,
\begin{equation}\label{equ:maincon2}
    c_2(X, X_{\tt des}, \theta) = U_{\tt max}^2 - \| U \|_2^2 \geq 0, \quad 
    U = -G^{\sf T}(X,\theta) P (X - X_{\tt des}),
\end{equation}
and ensures that the spacecraft relative acceleration due to thrust, $U$, remains below a specified value, $U_{\tt max}$, in magnitude. This constraint assumes that the spacecraft has a single orbital maneuvering thruster and that the attitude of the spacecraft is changed by an attitude control loop to accurately realize the commanded thrust direction.  
We also assume that the standard pulse width pulse frequency modulation \cite{wie1998space} is employed to realize continuous values of thrust magnitude if on-off thrusters are employed onboard of the spacecraft.  Furthermore,  we neglect the minimum impulse bit lower limit on positive thrust magnitude as below this limit the thrust is set to zero, while with our approach, based on GVEs and invariance, and given the drift-free form of (\ref{equ:main1}), the constraints will not be violated.  In fact, this is an additional reason for why the minimum distance constraint (\ref{equ:maincon1}) is imposed on the radius of periapsis.

The third constraint limits the eccentricity of the spacecraft orbit d above a specified minimum value, $e \geq e_{\tt min}$, where $e_{\tt min} \geq 0$:
\begin{equation}\label{equ:maincon3}
    c_3(X) = e-e_{\tt min} \geq  0.
\end{equation}
This constraint ensures that the eccentricity does not approach zero too closely where GVEs have a singularity, thereby preserving the validity of the model.  

Other constraints could be imposed depending on the mission objectives; such constraints can be similarly handled with the proposed reference governor approach.

\subsection{Multi-Step Incremental Reference Governor}\label{sec:refgov}

The reference governor\cite{garone2017reference} modifies the command $X_{\tt des}$ in transients to enforce the constraints.  That is, the reference command $X_{\tt des}$ in the control law (\ref{equ:main3}) is replaced by a modified reference command, $\Tilde{X}$, which does not cause constraint violations.  The control law (\ref{equ:main3}) thus becomes
\begin{equation} \label{equ:main3mod}
    U(t)=-G^{\sf T}(t)P(X(t)-\Tilde{X}(t)).
\end{equation}

The Incremental Reference Governor (IRG) \cite{tsourapas2009incremental} is a variant of the reference governor that tests the feasibility of a sufficiently small increment in $\Tilde{X}$ towards $X_{\tt des}$ at each discrete-time instant. If feasible, i.e., it is possible to infer that the constraints will be satisfied over a semi-infinite  prediction horizon with $\Tilde{X}$ being constantly applied, this increment is implemented. If not feasible, $\Tilde{X}$ is kept unchanged. In the previously considered fuel cell application of IRG\cite{tsourapas2009incremental}, the constraints were imposed to keep the system in the closed-loop region of attraction when the load was changing; the characterization of this region of attraction was computed numerically offline based on the reduced order model.

To develop the IRG for our orbital transfer application, we let $t_k$ denote a time instant at which IRG modifies the reference command and $X(t_k)$ denote the spacecraft state at that time instant.  
We re-write (\ref{equ:main2.5}) as
\begin{equation}\label{equ:main2.5.5}
V(X,\Tilde{X},P)=\frac{1}{2} (X-\Tilde{X})^{\sf T} P (X-\Tilde{X}),
\end{equation}
to make the dependence on the modified reference command, $\Tilde{X}$, and the weight matrix $P$ explicit.  The sub-level set of $V$ is defined as
\begin{equation}
Q\big(\Tilde{X},P,X(t_k)\big)=\bigg\{X \in \mathbb{R}^5:~V(X,\Tilde{X},P) \leq V(X(t_k),\Tilde{X},P)\bigg\}.
\end{equation}
Note that with the reference command maintained constant, $\Tilde{X}(t)=\Tilde{X}$ for $t \geq t_k$, and given that $\dot{V} \leq 0$ per (\ref{equ:main4.5}) with $X_{\tt des}=\Tilde{X}$, the closed-loop trajectory of (\ref{equ:main1}), (\ref{equ:main3mod}) will remain within 
the sublevel set for all future times, i.e.,
$$V(X(t),\Tilde{X},P) \leq V(X(t_k),\Tilde{X},P), \quad \mbox{ for $t \geq t_k$}.$$
If the sublevel set is constraint admissible (i.e., constraints are satisfied by all points of this set), then the trajectory is guaranteed to satisfy the constraints.

Testing if the invariant set is constraint admissible can be accomplished using optimization. Given the state at the time instant $t_k$, $X(t_k)$, define,
\begin{equation}\label{equ:mycon1}
    c_1^*(\Tilde{X},P,X(t_k)) = \min_{X\in Q(\Tilde{X},P,X(t_k))} c_1(X),
\end{equation}
\begin{equation}\label{equ:mycon2}
    c_2^*(\Tilde{X},P,X(t_k)) = \min_{X\in Q(\Tilde{X},P,X(t_k)),~0 \leq \theta < 2 \pi} c_2(X,\Tilde{X},\theta),
\end{equation}
and
\begin{equation}\label{equ:mycon3}
    c_3^*(\Tilde{X},P,X(t_k)) = \min_{X\in Q(\Tilde{X},P,X(t_k))} c_3(X).
\end{equation}
We refer to a reference command $\Tilde{X}$ as admissible at the time instant $t_k$ if $c_1^*(\Tilde{X},P,X(t_k)) \geq 0$, $c_2^*(\Tilde{X},P,X(t_k)) \geq 0$ and   $c_3^*(\Tilde{X},P,X(t_k)) \geq 0$.  
In our implementation of m$^2$IRG in this Section, given high dimensionality of $\tilde{X}$ and $X(t_k)$ and varying $P$, we resort to checking admissibility of $\Tilde{X}$ by solving the optimization problems (\ref{equ:mycon1})-(\ref{equ:mycon3}) online rather than pre-computing offline
and storing approximations of $c_1^*(\Tilde{X},P,X(t_k))$, $c_2^*(\Tilde{X},P,X(t_k))$ and  $c_3^*(\Tilde{X},P,X(t_k))$.
We will also  consider the implementation of m$^2$IRG using online predictions without solving (\ref{equ:mycon1})-(\ref{equ:mycon3}) online in Section~\ref{sec:8}.

At each time instant $t_k$, $k>1$, the multi-step IRG evaluates the admissibility of $n>0$ reference commands, $\Tilde{X}^j(t_k)$ $j=1,\cdots,n$, generated by the following rule,
\begin{align}
& \Tilde{X}^0(t_k)=\Tilde{X}(t_{k-1}) \nonumber \\
& \Tilde{X}^{j}(t_k)=\Tilde{X}^{j-1}(t_k)+\Delta_k^j E_k^j (X_{\tt des} - \Tilde{X}^{j}(t_k)),\quad j>0,
\end{align}
where $\Delta_k^j>0$ is a scalar step-size and $E_k^j \in \mathbb{R}^{5 \times 5}$ determines the direction of adjustment. 

Note that if an admissible $\Tilde{X}(0)$ exists at the initial time (typically the choice $\Tilde{X}(0)=X(0)$ suffices), then the constraints are guaranteed to be satisfied for all future times. Based on similar analysis as in the exactly implemented command governor case \cite{garone2022command}, it can also be shown that the choice $E_k^j={\tt diag}([1,1,1,1,1,1])$, $\Tilde{X}(0)=X(0)$ and {\em strict} constraint admissibility of all points along the line segment $\lambda \Tilde{X}(0)+(1-\lambda) X_{\tt des},~0 \leq \lambda \leq 1$, are sufficient to ensure the convergence of $\Tilde{X}(t_k)$ to $X_{\tt des}$ in a finite number of steps, provided the steps $\Delta_k^j>0$ are chosen to be sufficient small and the persistence of excitation condition holds, which ensures closed-loop convergence to constant reference commands. 

In practice, we adopt a more flexible version of the algorithm to mitigate situations in which the above ``straight line'' descent 
results in slow convergence.  In this case,  $E_k^j$ is  chosen as
$$ E_k^j = \left\{ \begin{array}{ll} {\tt diag} (e_l) & \mbox{if }
k ~{\tt mod}~6=l,~l=0,\cdots,4, \\
{\tt diag}([1,1,1,1,1]) & \mbox{if } 
k ~{\tt mod}~6=5,
\end{array}
\right.
$$
where $e_l$ is the $l$th unit vector in $\mathbb{R}^5$, $l=1,\cdots,5$.
This approach uses coordinate descent interlaced with periodic ``straight line'' descent to $X_{\tt des}$ which remains feasible given the geometry of the constraints in our problem, provided the reference commands satisfying $c_1(\Tilde{X}^j(t_k))<\epsilon$ and $c_3(\Tilde{X}^j(t_k))<\epsilon$ for a sufficiently small $\epsilon>0$ (i.e., too close to the boundary of the constrained region) are declared inadmissible.

The multi-step IRG algorithm functions as follows. Initially, a nominal step size is chosen, $\Delta_k^1=\Delta$. If $\Tilde{X}^1(t_k)$ is inadmissible then the step size is reduced, $\Delta_k^2=\gamma \Delta_k^1$, $0<\gamma<1$.
If $\Tilde{X}^2(t_k)$ is inadmissible then $\Tilde{X}(t_k)=\Tilde{X}(t_{k-1})$ is passed to the Lyapunov controller and the algorithm stops.  If $\Tilde{X}^{j-1}(t_k)$ is admissible, we set $\Delta_k^j=\Delta_{k}^{j-1}$.  If $\Tilde{X}^{j}(t_k)$ is inadmissible then $\Tilde{X}(t_k)=\Tilde{X}^{j-1}(t_k)$ and the algorithm stops.  

Note that the IRG algorithm can gracefully handle convergence failures of the numerical optimizers used to solve the low dimensional nonlinear optimization problems (\ref{equ:mycon1})-(\ref{equ:mycon3}).  In such cases, the choice $\Tilde{X}^j(t_k)=\Tilde{X}^{j-1}(t_k)$ remains feasible.

\subsection{Multi-Step Multi-Mode Incremental Reference Governor}

The m$^2$IRG
represents an extension of Multi-Step IRG to allow the choice of $P$ matrix in the Lyapunov function (\ref{equ:main2.5.5}) 
and in the Lyapunov controller (\ref{equ:main3}) to vary.
Specifically, assuming $P \in \mathcal{P}=\{P^1,\cdots,P^m\}$,
where $P^l=(P^l)^{\rm T} \succ 0$, $l=1,\cdots,m,$ the criterion for
admissibility of a reference command $\Tilde{X}^j$ can be modified to existence of $P \in \mathcal{P}$ such that 
$$c_j^*(\Tilde{X}^j,P,X(t_k)) \geq 0, \quad j=1,2,3.$$  
The use of $\mathcal{P}$, however, increases the number of optimization problems which need to be solved onboard.

A different strategy that does not require as much computations is motivated by the observation that the primary constraints are (\ref{equ:maincon1}) and (\ref{equ:maincon3}).  This strategy involves identifying  $P \in \mathcal{P}$ that best matches the current values of the semi-major axis $a(t_k)$ and eccentricity $e(t_k)$.  The basic idea is to choose $P \in \mathcal{P}$ so that the invariant set, the shape and orientation of which are determined by $P$, is stretched along the region allowed by the constraints so that to facilitate rapid progress of $\Tilde{X}$ towards $X_{\tt des}$. This best $P$ is denoted by $P_{\tt des}(t_k)$.  The IRG algorithm can then be modified so that its first step is to test whether $P_{\tt des}(t_k)$ and $\Tilde{X}(t_{k-1})$ are admissible  at the time instant $t_k$. If admissible, then $P(t_k)=P_{\tt des}(t_k)$ and the algorithm continues to search for $\Tilde{X}(t_k)$. If not admissible, then $P(t_k)=P(t_{k-1})$ and the algorithm exits with $\Tilde{X}(t_{k})=\Tilde{X}(t_{k-1}).$  Under the persistence of excitation conditions, $X(t) \to \Tilde{X}(t)$ if $\Tilde{X}$ and $P$ do not vary; hence eventually, the switch to $P_{\tt des}(t_k)$ becomes possible. In the actual implementation, we found that continuing to search for an improvement in $\Tilde{X}(t_{k})$ even when $P(t_k) $ cannot be updated to $P_{\tt des}(t_k)$ resulted in faster response; hence we adopted this latter choice.

\section{Numerical Results with Lyapunov Function Based Incremental Reference Governor}\label{sec:5}
The   m$^2$IRG is first tested in orbital transfer maneuvers around Earth where the constraints
(\ref{equ:maincon1})-(\ref{equ:maincon3})
are defined with
$r_{\tt min}=6628$, $U_{\tt max}=1.25\times 10^{-3} \frac{\text{km}}{\text{s}^2}$, $e_{\tt min}=10^{-6}$. Lower acceleration limit values could be handled similarly, however, the maneuver time and simulation time will increase.

The m$^2$IRG implementation uses $\Delta=0.01$, $\gamma=\frac{1}{5}$, $n=12$.
The optimization problems (\ref{equ:mycon1})-(\ref{equ:mycon3}) were solved using the function {\tt fmincon} of Matlab with the initial guess set to $X(t_k)$, $\theta(t_k)$ and maximum number of iterations set to $600$.
Additional acceptance logic for the computed $\tilde{X}^i(t_k)$ has been added to handle occasional {\tt fmincon} numerical convergence issues. This logic checks that {\tt fmincon} has (approximately) converged and that the constraint (\ref{equ:maincon2})
at the current time instant $t_k$ is satisfied.  
The reference command is updated by m$^2$IRG every $t_{k+1}-t_k=15$ min.

A set of matrices $\mathcal{P}=\{P^1,P^2,P^3\}$ has been constructed  as follows.
First, consistently with our Lyapunov controller calibration guidelines in Section~\ref{sec:Lyapunov cont}, we manually tuned the diagonal entries of a single diagonal $P$ matrix to improve the speed of transient response of individual orbital elements (larger diagonal element speeds up the transient response of the corresponding orbital element) though simulations. This led to the matrix
$$P^0={\tt diag}\left( [\begin{array}{ccccc} 5 \times 10^{-11}, & 0.1, & 5 \times 10^{-3}, & 7.5 \times 10^{-3}, & 5 \times 10^{-4}  \end{array}] \right).$$
Note that the magnitudes of these diagonal elements are quite different, but these differences are primarily attributed to differences in the ranges of the orbital elements involved (e.g., $a$ is in 1000's of km, while $e$ is less than $1$).

\begin{figure}[h]
        \centering
        \includegraphics[width = 8 cm]{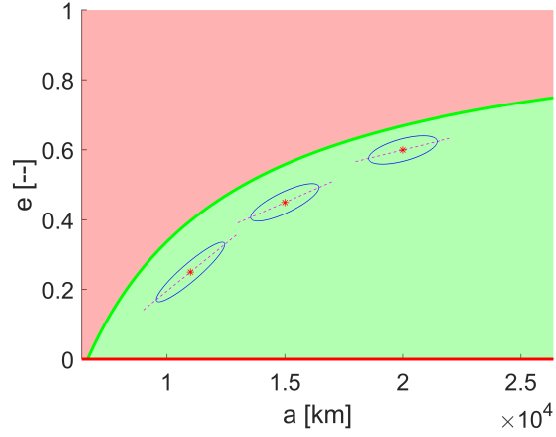}
        \caption{Illustration of cross sections of sublevel sets of the Lyapunov function corresponding to $P^1$, $P^2$ and $P^3$.}
        \label{fig:invariantset1}
    \end{figure}

We then considered rotations of the cross sections of the sublevel sets of the Lyapunov function
by the $a$-$e$ plane to orient them along the boundary prescribed by the constraint (\ref{equ:maincon1}) at three different values of $a$, $a_1=2 \times 10^4$, $a_2=1.5 \times 10^{4}$ and $a_3=1.1 \times 10^4$.  See 
Figure~\ref{fig:invariantset1}.
This reorientation of the cross sections facilitates larger
changes in $\tilde{X}(t_k)$.
By characterizing the tangent to the constraint boundary (\ref{equ:maincon1}), the angle by which the cross section needs to be rotated is 
$\alpha_j=\arctan(\frac{r_{\tt min}}{a_j^2})$, $j=1,2,3$.
If $P^j_{(2,2)}$ denotes the upper $2 \times 2$ block of the $P^j$ matrix, $j=1,2,3$ then
$$P_{(2,2)}^j=R(\alpha_j) P^0 R^{\sf T}(\alpha_j),\quad 
R(\alpha)=\left[\begin{array}{cc} \cos(\alpha) & \sin(\alpha) \\ - 
\sin(\alpha) &  \cos(\alpha) \end{array} \right].$$
With this process we generated three matrices,$P^1$,$P^2$,$P^3$:
 $$P^1=\left[\begin{array}{ccccc} 
 7.7456 \times 10^{-11} & -1.656999999 \times 10^{-6} & 0 & 0 &0 \\
  -1.656999999 \times 10^{-6} &  0.099999999972544 & 0 & 0 \\
  0 & 0 & 5 \times 10^{-3} & 0 & 0 \\
  0& 0& 0 & 7.5 \times 10^{-2} & 0 \\
  0 & 0 &0 &0 & 5 \times 10^{-4}
\end{array} \right],
$$

$$P^2=\left[\begin{array}{ccccc} 
 2.61856 \times 10^{-10} & -4.602777768 \times 10^{-6} & 0 & 0 &0 \\
  -4.602777768 \times 10^{-6}  & 0.099999999788144 &0 & 0 & 0 \\
  0 & 0 & 5 \times 10^{-3} & 0 & 0 \\
  0& 0& 0 & 7.5 \times 10^{-2} & 0 \\
  0 & 0 &0 &0 & 5 \times 10^{-4}
\end{array} \right],
$$

$$P^3=\left[\begin{array}{ccccc} 
  1.066157 \times 10^{-9} & -1.0080463648 \times 10^{-5} & 0 & 0 &0  \\
 -1.0080463648 \times 10^{-5}  & 0.099999998983843 & 0 & 0 &0 \\
  0 & 0 & 5 \times 10^{-3} & 0 & 0 \\
  0& 0& 0 & 7.5 \times 10^{-2} & 0 \\
  0 & 0 &0 &0 & 5 \times 10^{-4}
\end{array} \right].
$$

The logic for selecting $P_{\tt des}(t_k) \in \mathcal{P}$ was as follows:
If $1.5 \times 10^4 \leq a(t_k)$ then
$P_{\tt des}(t_k)=P^1$.  If $1.1 \times 10^4 \leq a(t_k) < 1.5 \times 10^4$ then
$P_{\tt des}(t_k)=P^2$.  If $a(t_k) < 1.1 \times 10^4$
then
$P_{\tt des}(t_k)=P^3$.

\subsection{Orbital transfer from higher orbit to lower orbit}
\label{sec:goingdown1}

\begin{figure}[h]
\centering
\subfloat[]{
    \includegraphics[width=4.5cm]{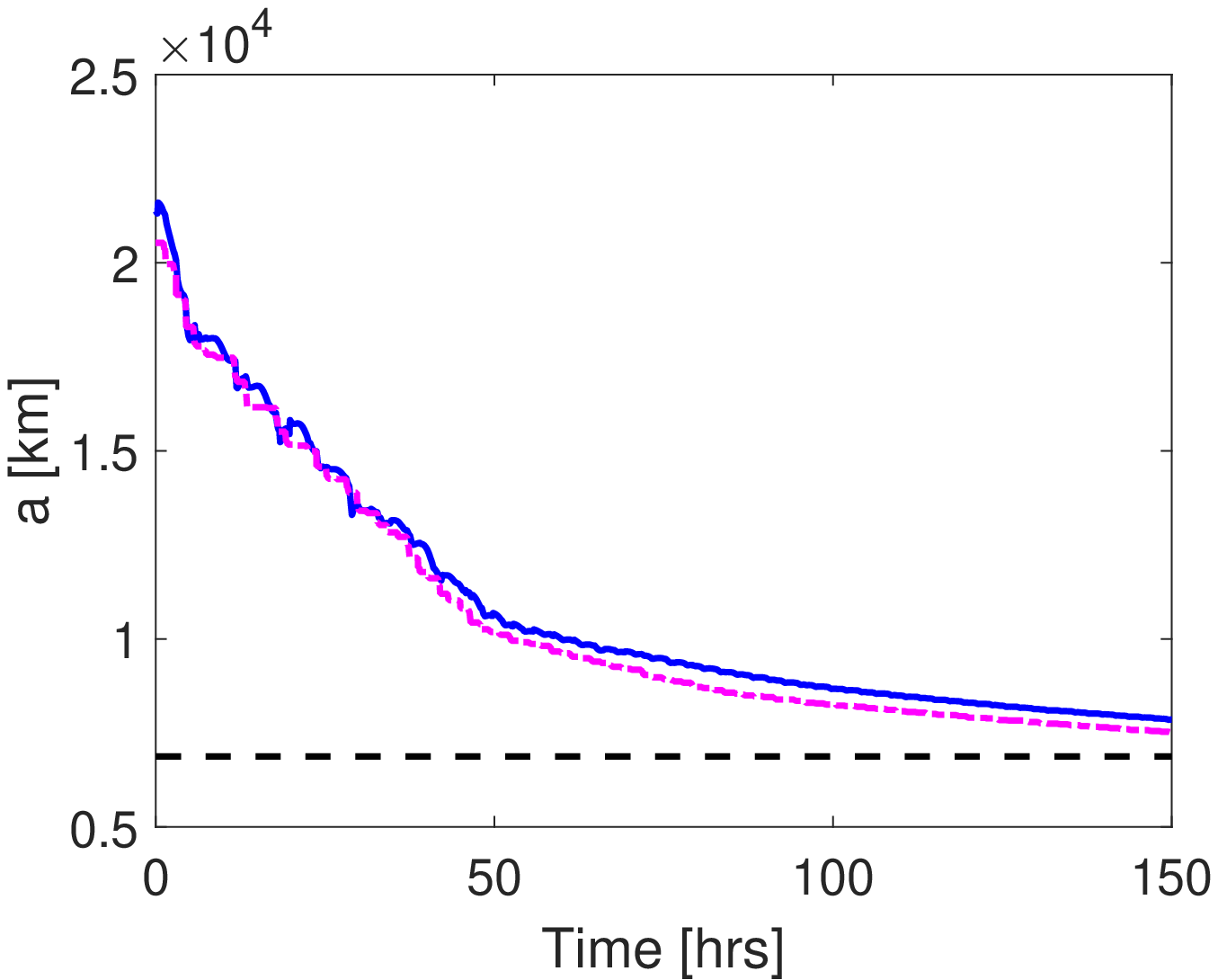}
}\qquad
\subfloat[]{
    \includegraphics[width=4.5cm]{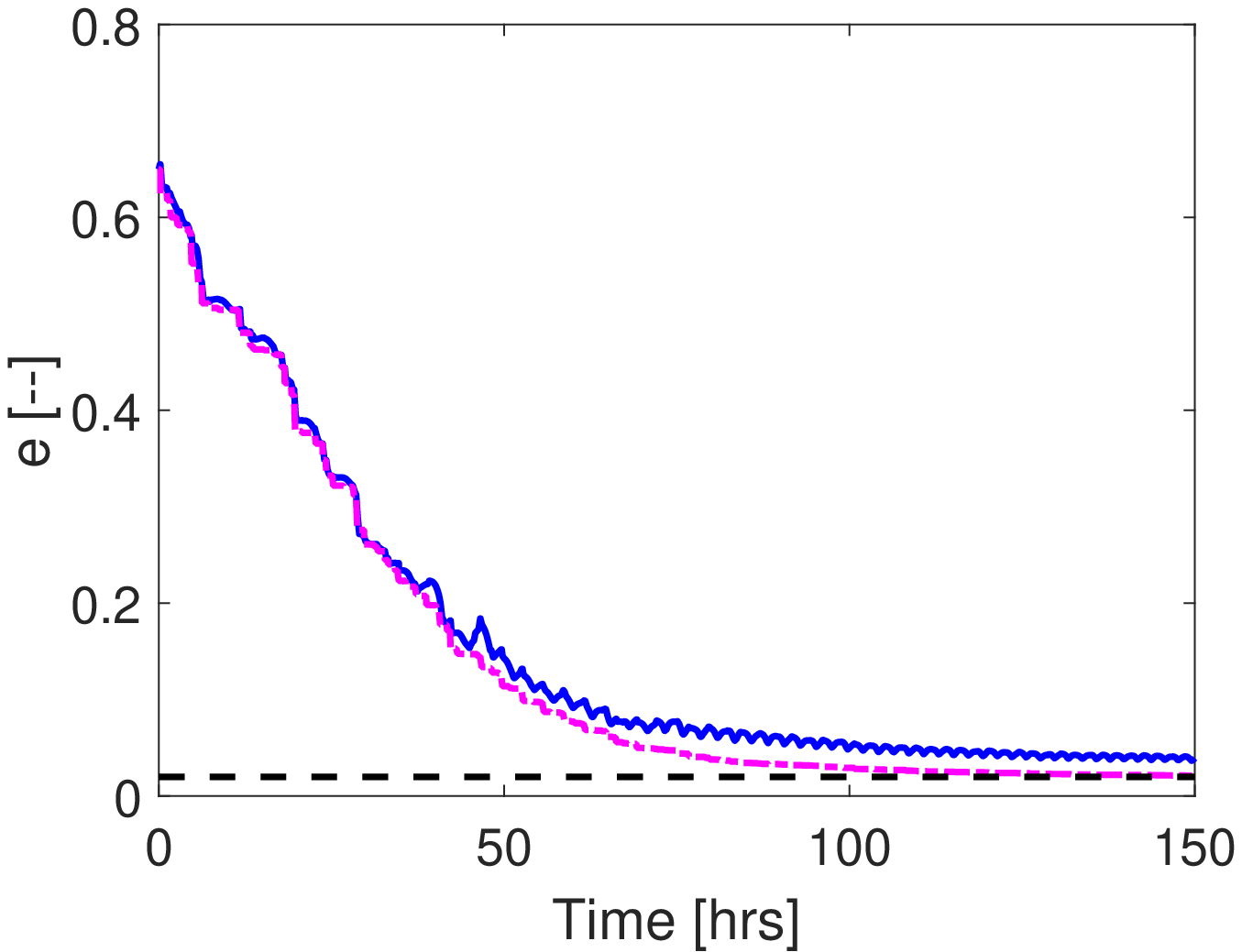}
\label{fig:callout0}}\qquad
\subfloat[]{
    \includegraphics[width=4.5cm]{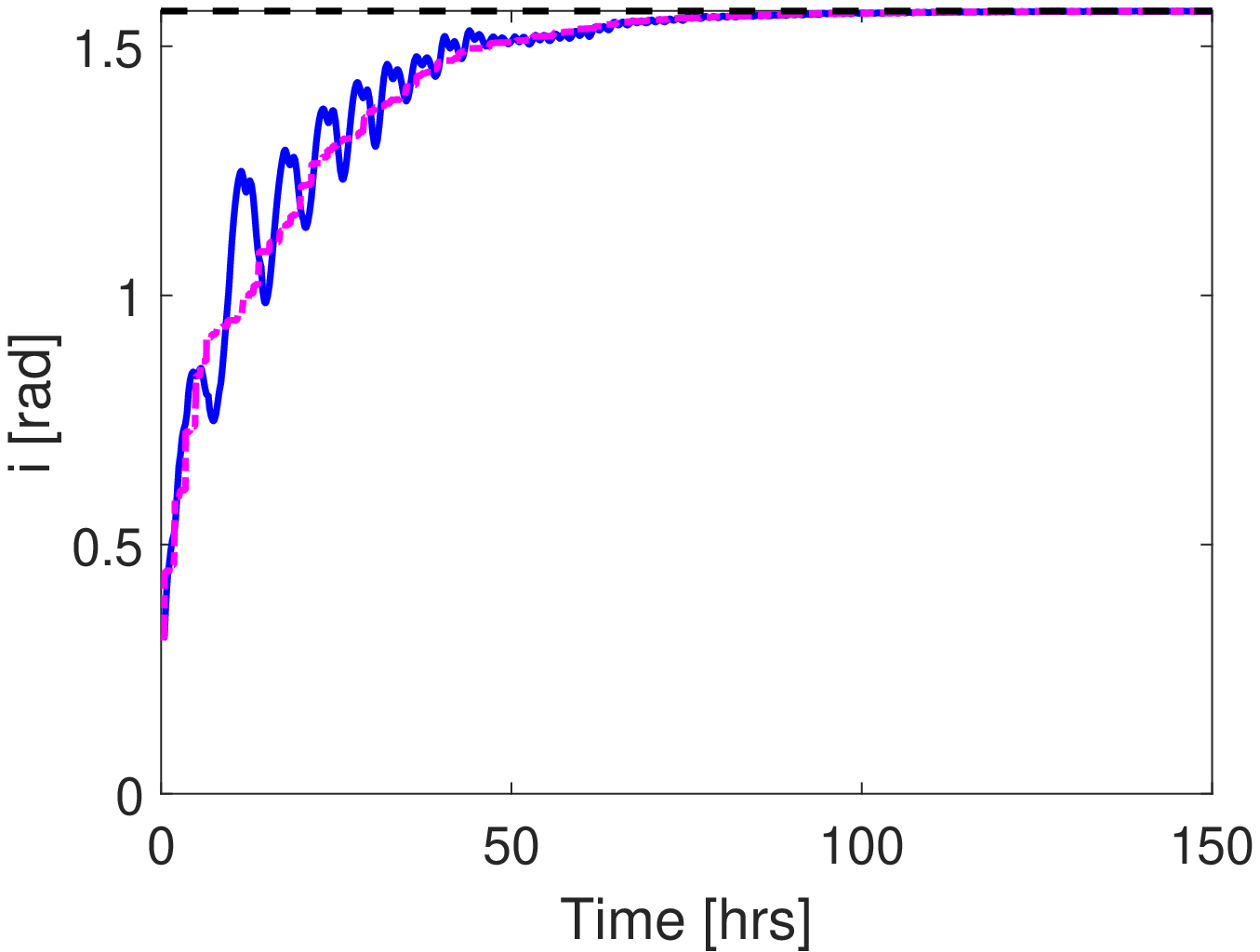}
}\vskip\baselineskip
\subfloat[]{
    \includegraphics[width=4.5cm]{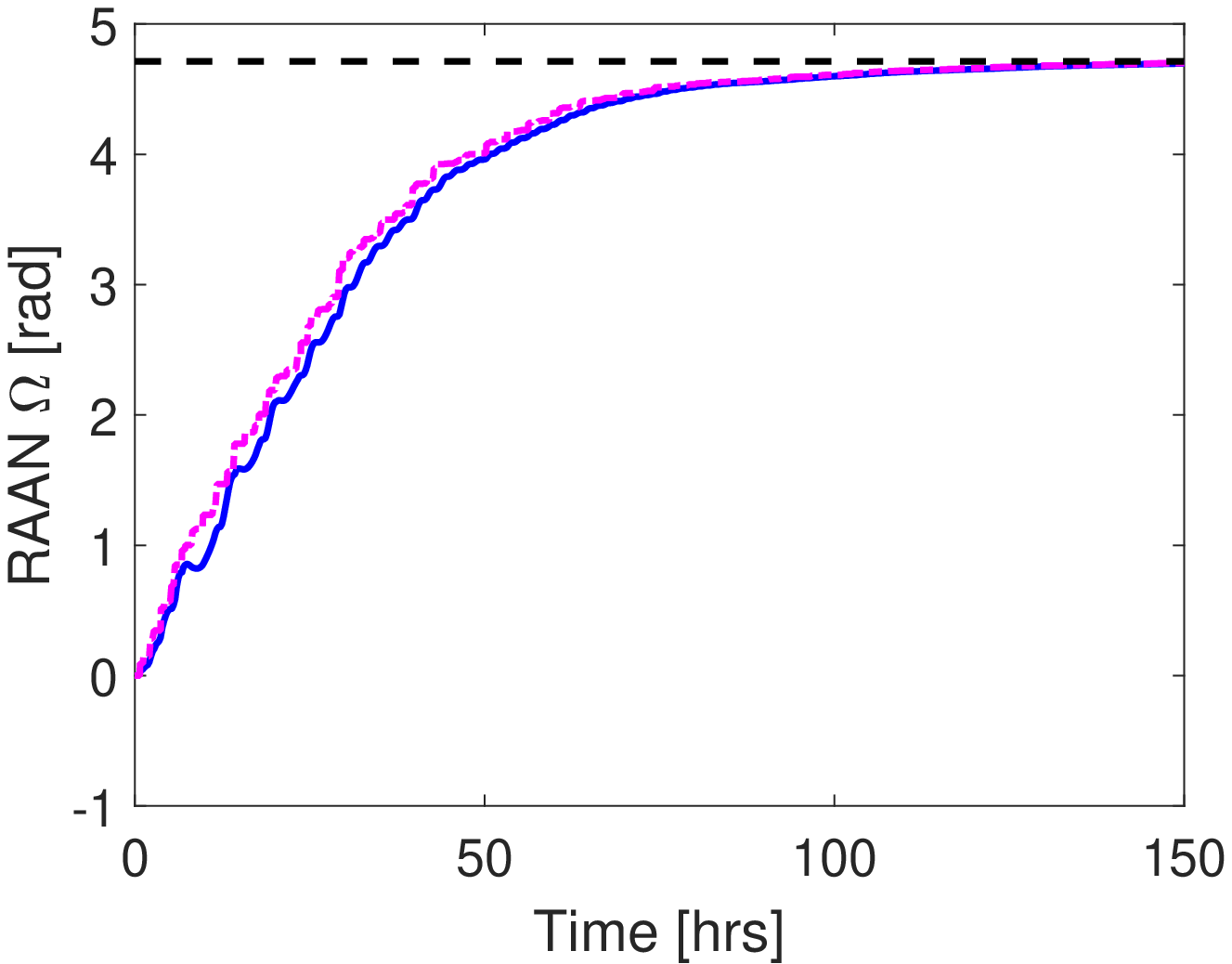}
}\qquad
\subfloat[]{
    \includegraphics[width=4.5cm]{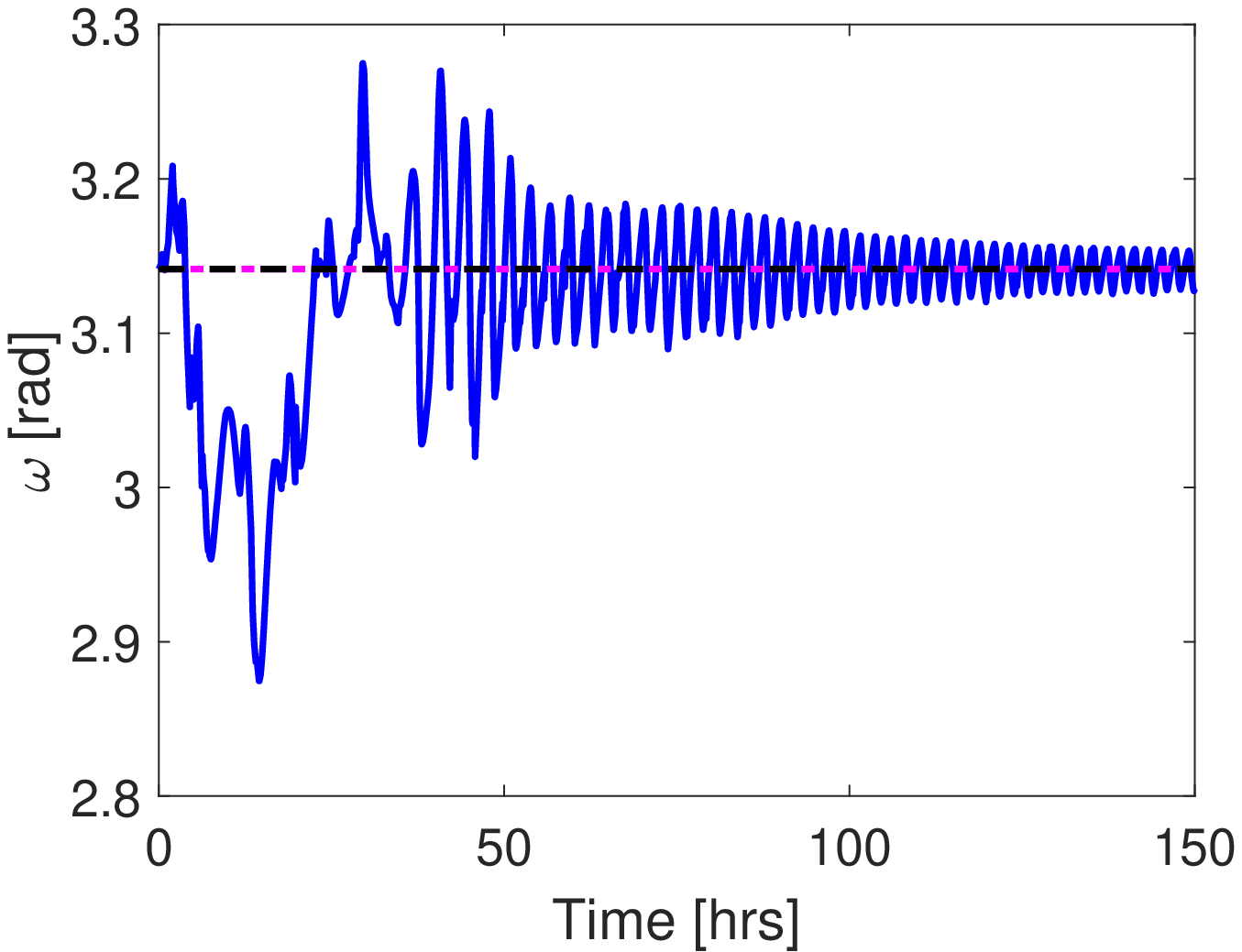}
}\qquad
\subfloat[]{
    \includegraphics[width=4.5cm]{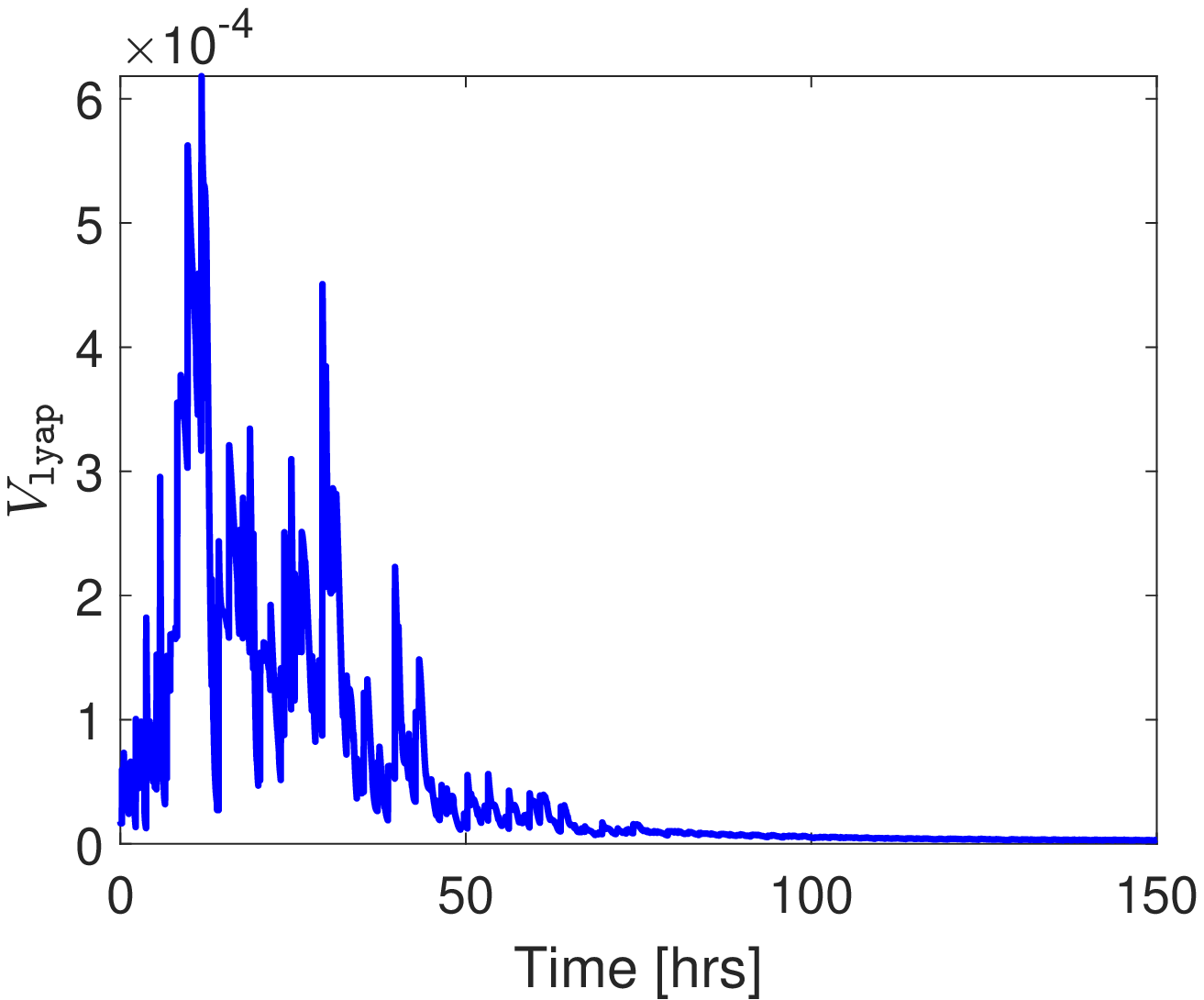}
    \label{subfig:Going down no barrier lyap}
}

\caption{Orbital transfer from a higher orbit to a lower orbit. The time histories of orbital elements are shown in blue, compared  with the final desired state (dashed, black) and the reference generated by the m$^2$IRG (magenta, dash-dotted). Figure~\ref{subfig:Going down no barrier lyap} is the time history of the Lyapunov function.}\label{fig:Going down no barrier}
\end{figure}

\begin{figure}[h]
\centering
\subfloat[]{
    \includegraphics[width=6cm]{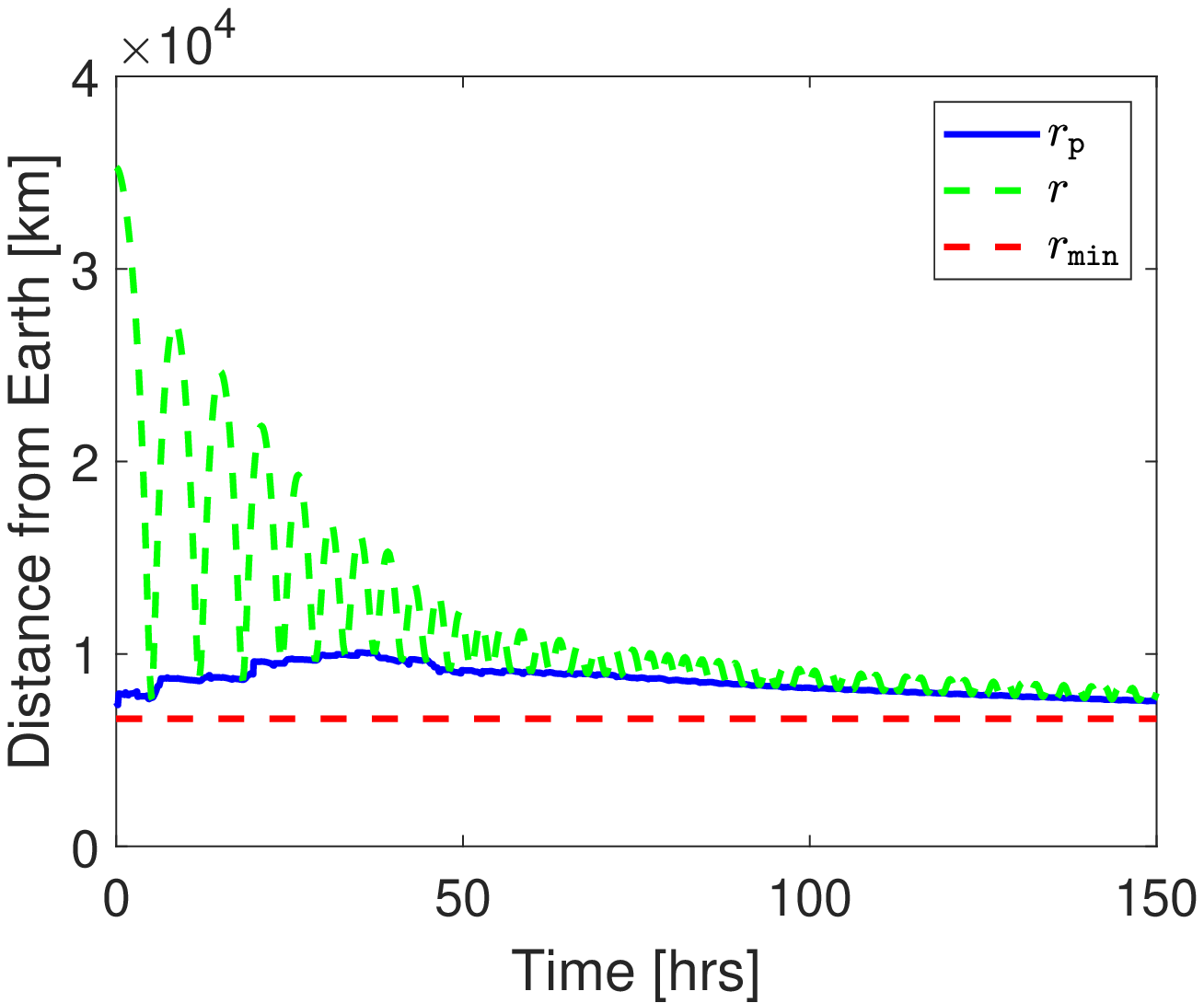}
}\qquad
\subfloat[]{
    \includegraphics[width=6cm]{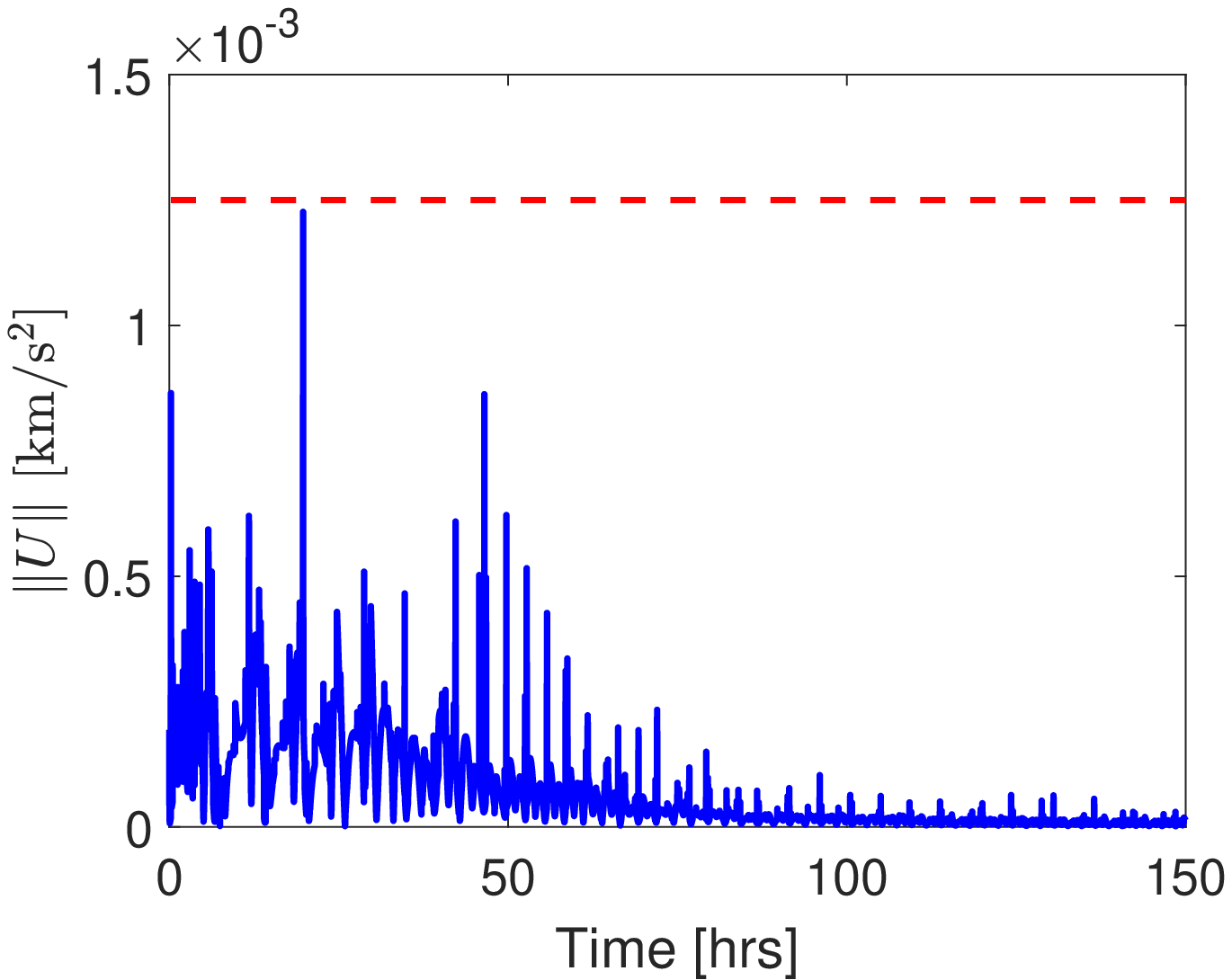}
}
\caption{Orbital transfer from a higher orbit to a lower orbit.  Left: The time histories of $r_p$, $r$ and $r_{\tt min}$ showing that the constraint (\ref{equ:maincon1}) is enforced. Right: The time histories of  $\|U\|$ (solid, blue) and of $U_{\tt max}$ (dashed, red)
showing  the constraint (\ref{equ:maincon2}) is enforced.}\label{fig:Going down coss}
\end{figure}

The maneuver from a higher orbit to a lower  orbit was simulated corresponding to the following initial state and desired state,
\begin{equation}
    X(0) = \left[21378, \ 0.65, \ \frac{\pi}{10}, \ 0, \ \pi, \ \pi \right]^{\sf T},
\end{equation}
to the target
\begin{equation}
    X_{\tt des} = [6878, \ 0.02, \ \frac{\pi}{2}, \ \frac{3\pi}{2}, \ \pi, \ 0]^{\sf T}.
\end{equation}
The time histories of the orbital elements, modified reference command, and of the Lyapunov function are shown in Figure~\ref{fig:Going down no barrier}. The target orbit is successfully transferred to. Figure~\ref{fig:Going down coss} shows that the constraints (\ref{equ:maincon1}) and (\ref{equ:maincon2}) are enforced during the maneuver.  From Figure~\ref{fig:callout0}), the constraint (\ref{equ:maincon3}) is enforced as well.

\subsection{Orbital transfer from lower orbit to higher  orbit }

With the same algorithm parameters, the maneuver from a lower orbit to a higher  orbit was simulated corresponding to the initial state and desired state that were interchanged as compared to the case in Subsection~\ref{sec:goingdown1}, i.e.,
\begin{equation}
    X(0) = [6878, \ 0.02, \ \frac{\pi}{2}, \ \frac{3\pi}{2}, \ \pi, \ 0]^{\sf T},
\end{equation}
\begin{equation}
    X_{\tt des} = \left[21378, \ 0.65, \ \frac{\pi}{10}, \ 0, \ \pi, \ \pi \right]^{\sf T}.
\end{equation}

\begin{figure}[h]
\centering
\subfloat[]{
    \includegraphics[width=4.5cm]{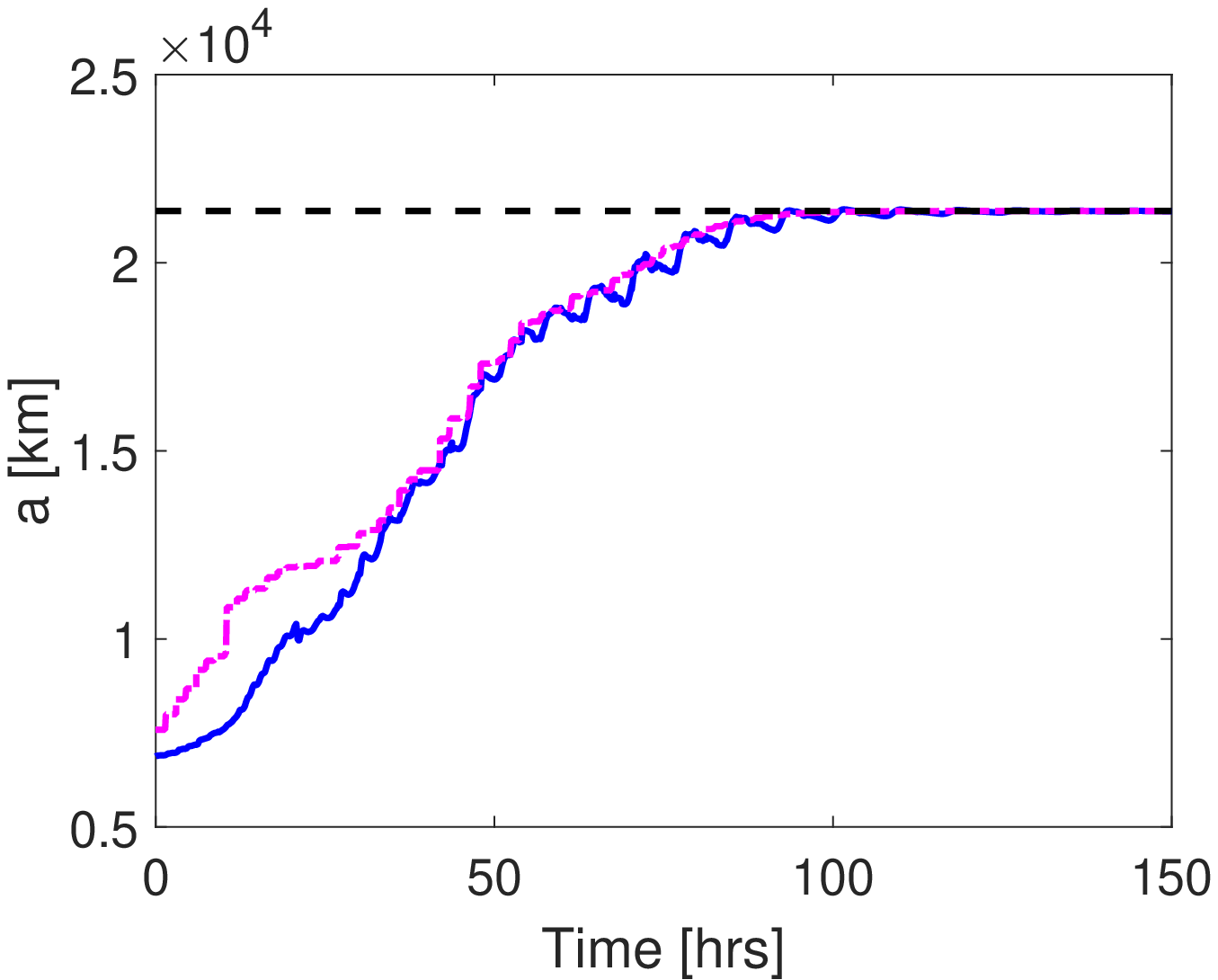}
}\qquad
\subfloat[]{
    \includegraphics[width=4.5cm]{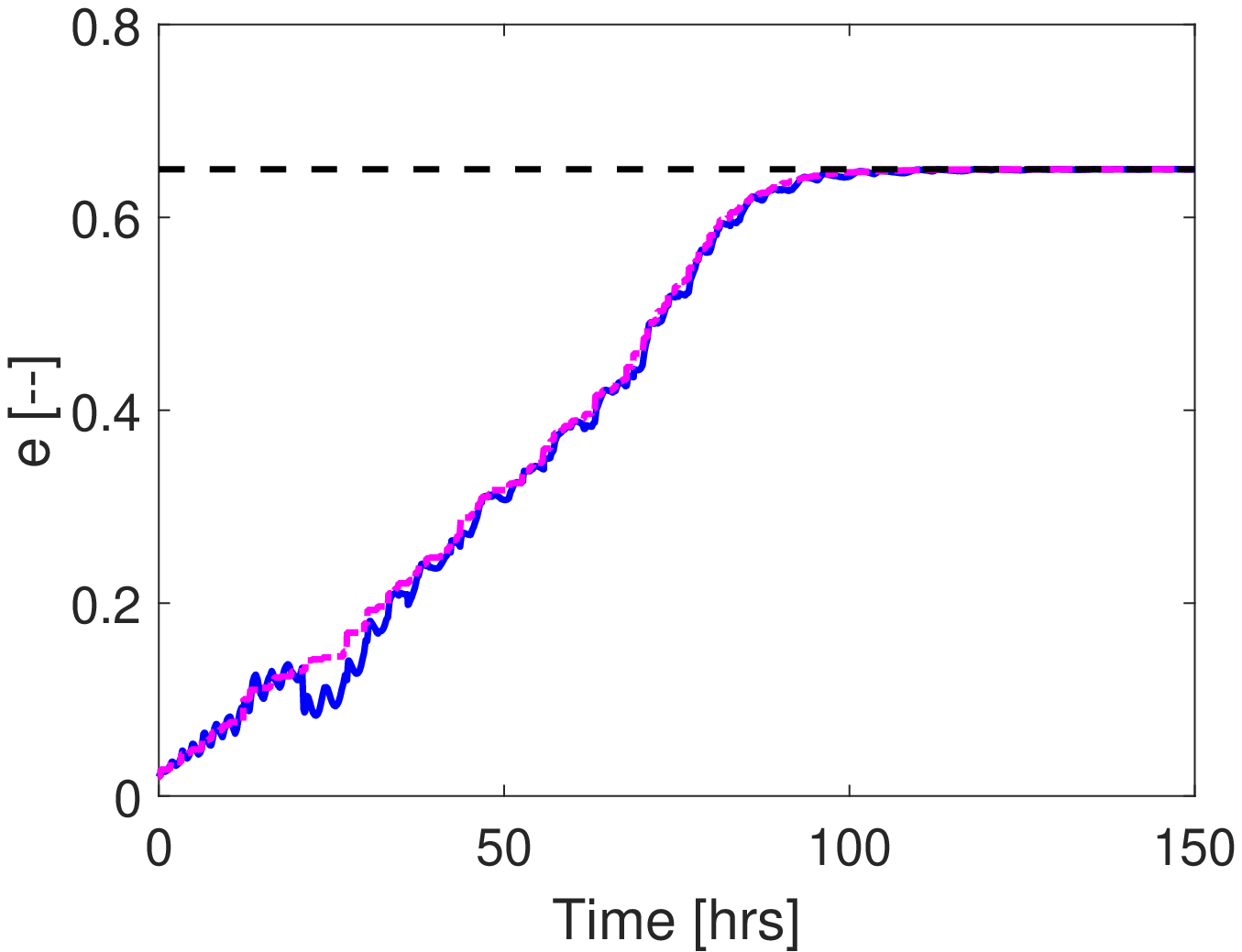}
\label{fig:callout1}}\qquad
\subfloat[]{
    \includegraphics[width=4.5cm]{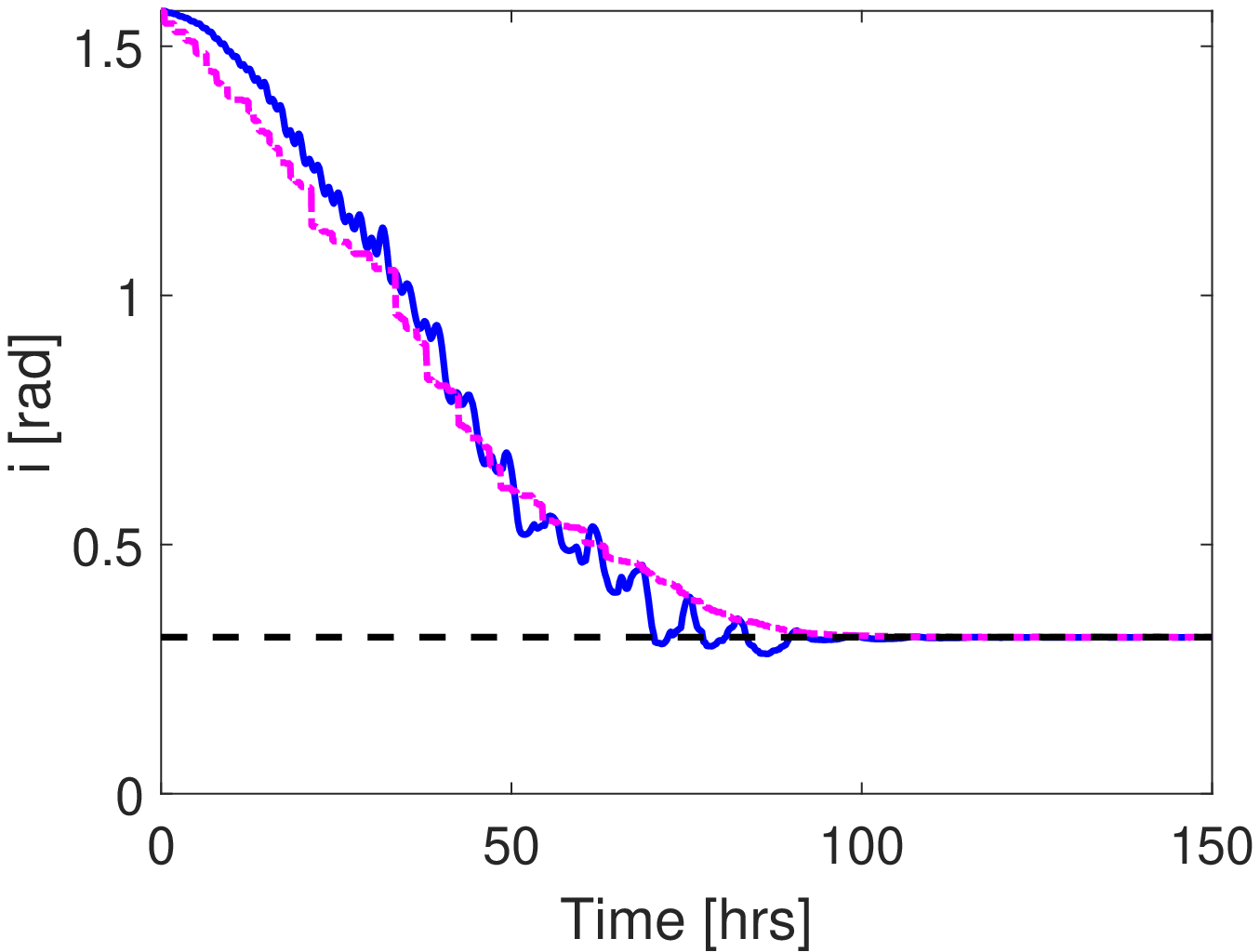}
}\vskip\baselineskip
\subfloat[]{
    \includegraphics[width=4.5cm]{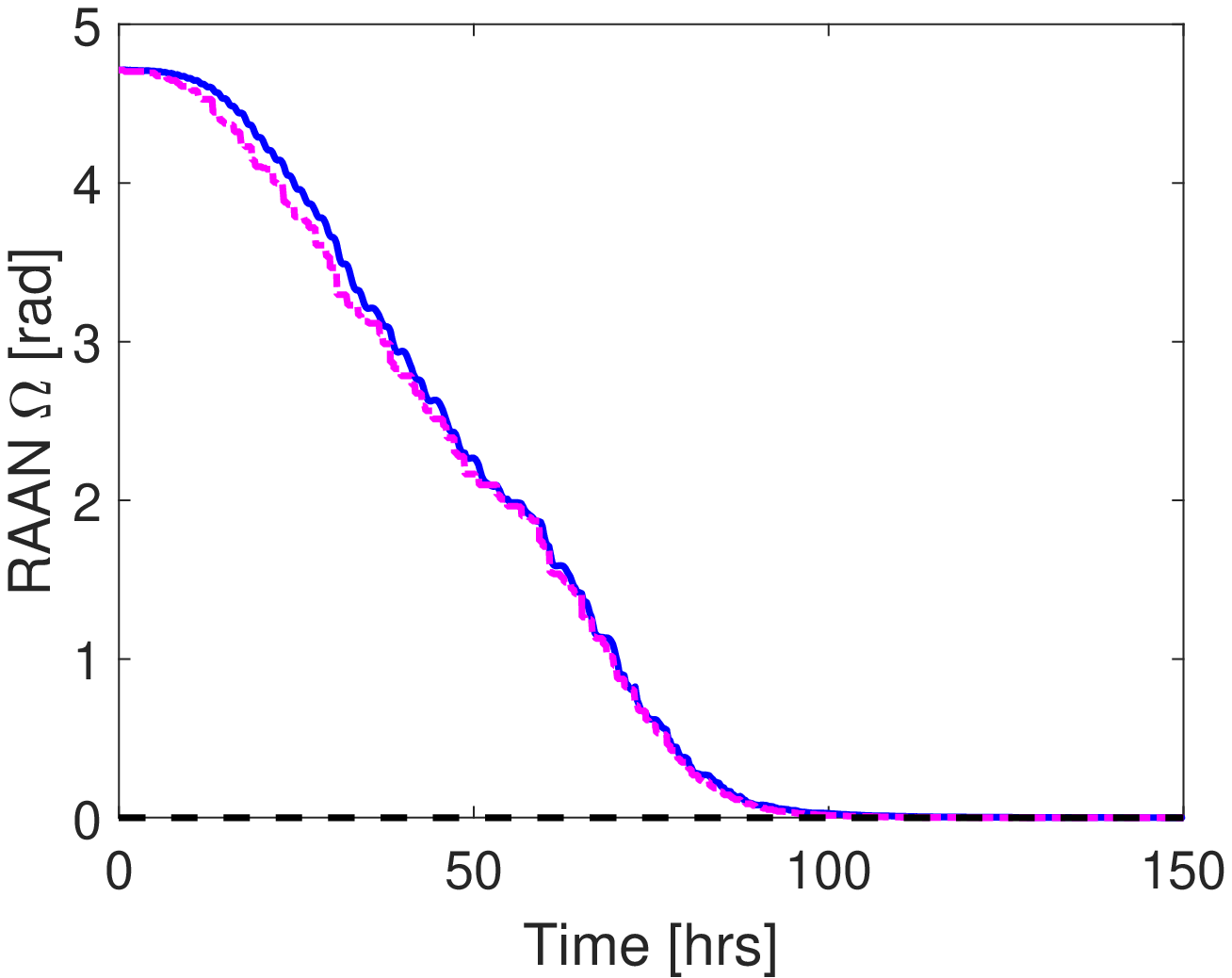}
}\qquad
\subfloat[]{
    \includegraphics[width=4.5cm]{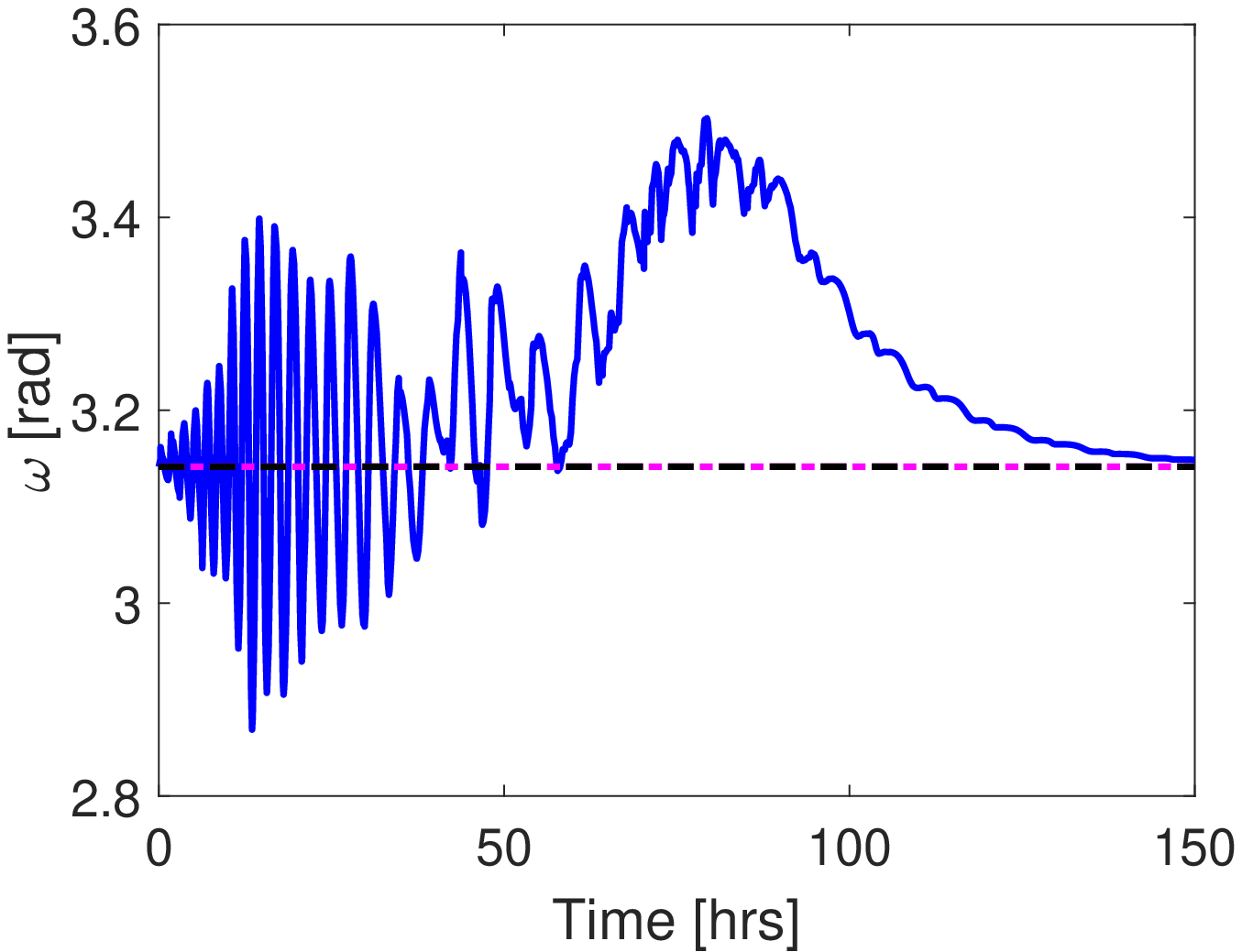}
}\qquad
\subfloat[]{
    \includegraphics[width=4.5cm]{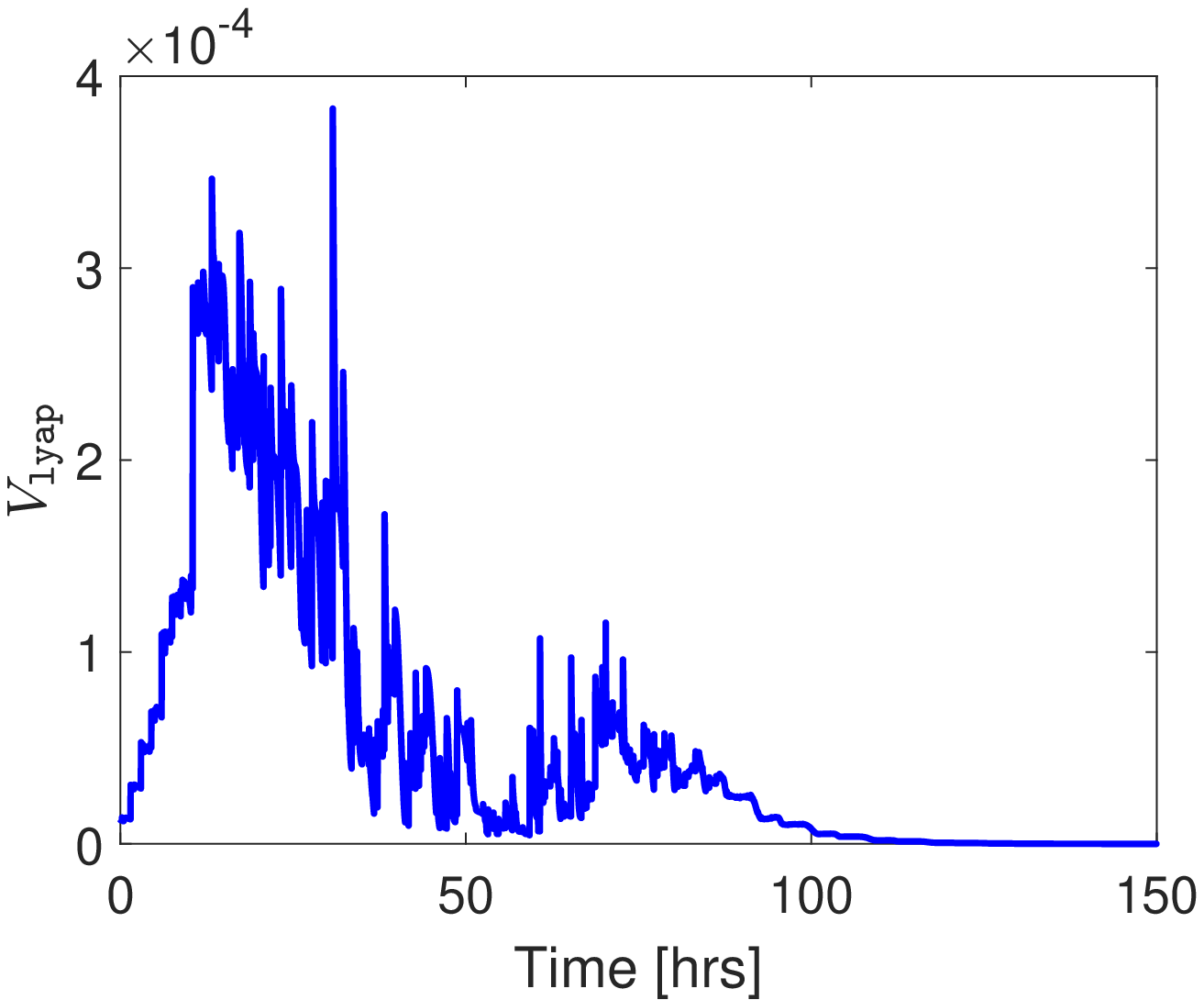}
    \label{subfig:Going up no barrier lyap11}
}
\caption{Orbital transfer from  a lower orbit to a higher orbit. The time histories of orbital elements are shown in blue, compared  with the final desired state (dashed, black) and the reference generated by the m$^2$IRG (dash-dotted, magenta). Figure~\ref{subfig:Going up no barrier lyap11} is the time history of the Lyapunov function.}\label{subfig:Going up no barrier lyap}
\end{figure}

\begin{figure}[h]
\centering
\subfloat[]{
    \includegraphics[width=6cm]{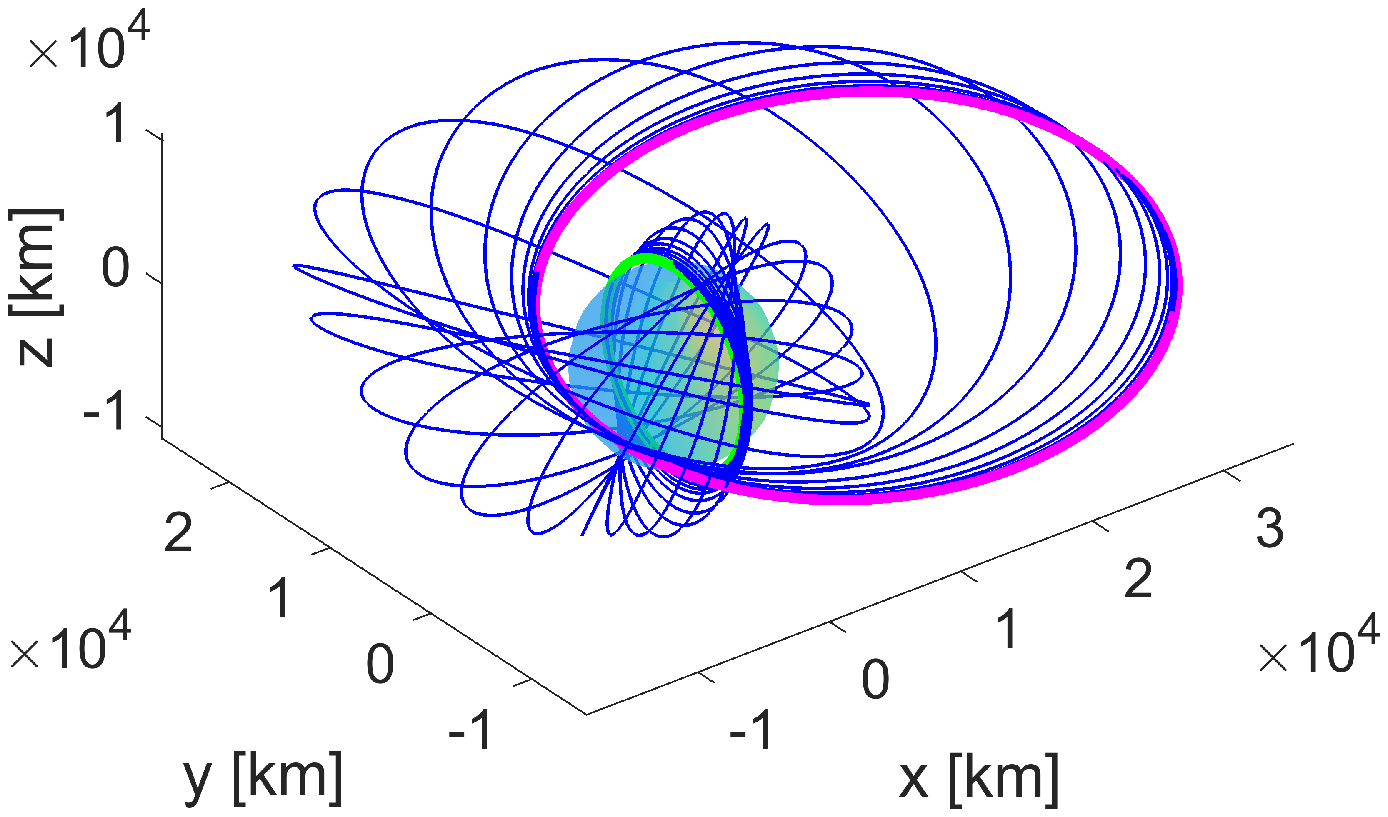}
}\qquad
\subfloat[]{
    \includegraphics[width=6cm]{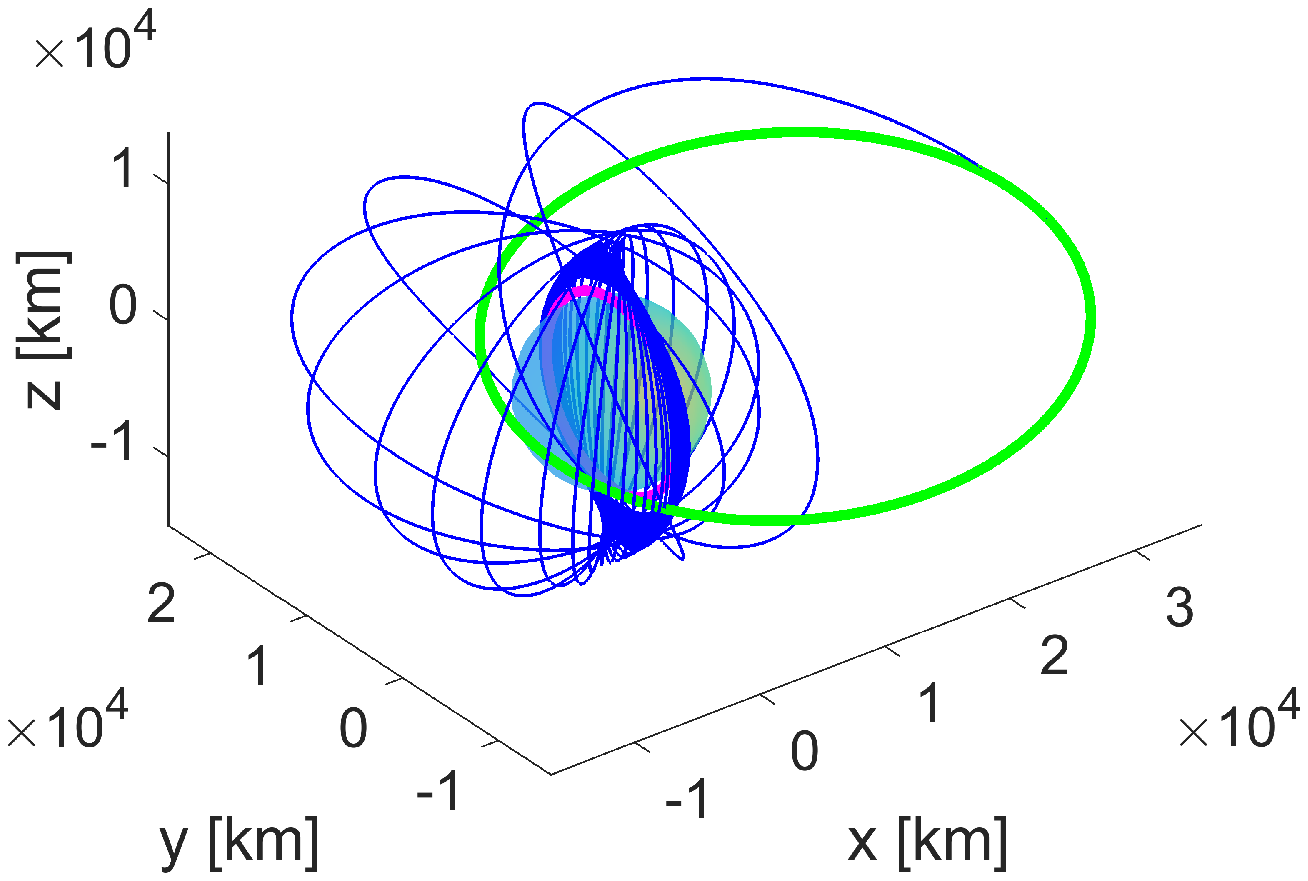}
}
\caption{Comparing three dimensional orbital transfer trajectories between a lower and a higher orbit. Initial orbit is shown by green, target orbit by magenta.}\label{fig:Traj 3d}
\end{figure}

The time histories of the orbital elements and of the Lyapunov function are shown in Figure~\ref{subfig:Going up no barrier lyap}.  The target orbit is successfully transferred to. Figure~\ref{fig:Going down coss} shows that the constraints (\ref{equ:maincon1}) and (\ref{equ:maincon2}) are enforced during the maneuver. From Figure~\ref{fig:callout1}, the constraint (\ref{equ:maincon3}) is enforced as well. Figure~\ref{fig:Traj 3d} compares three dimensional orbital transfer trajectories from lower to higher orbit and vice versa. It can be observed that m$^2$IRG is more cautious in making adjustments to reference commands when closer to Earth.

\begin{figure}[h]
\centering
\subfloat[]{
    \includegraphics[width=6cm]{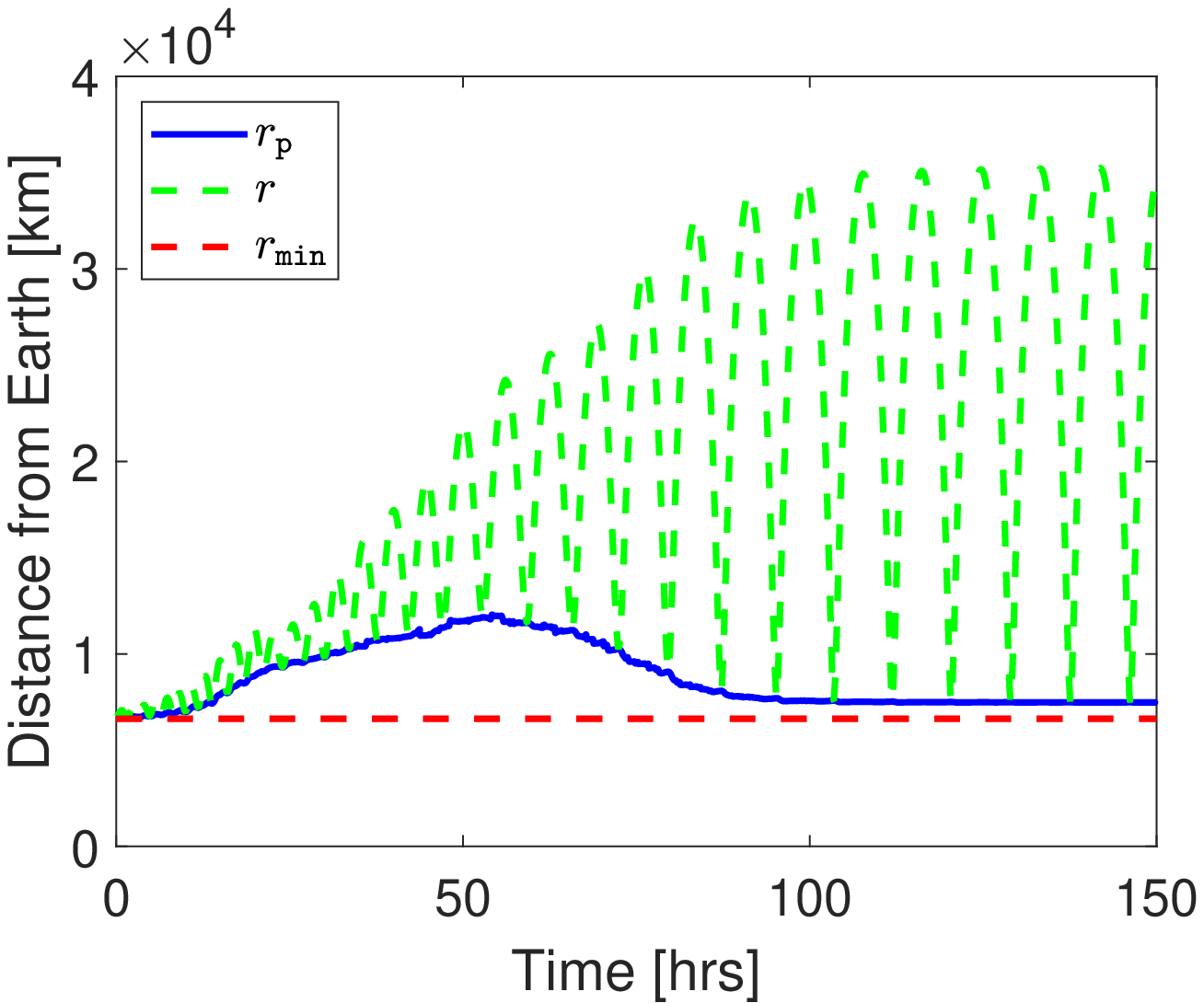}
}\qquad
\subfloat[]{
    \includegraphics[width=6cm]{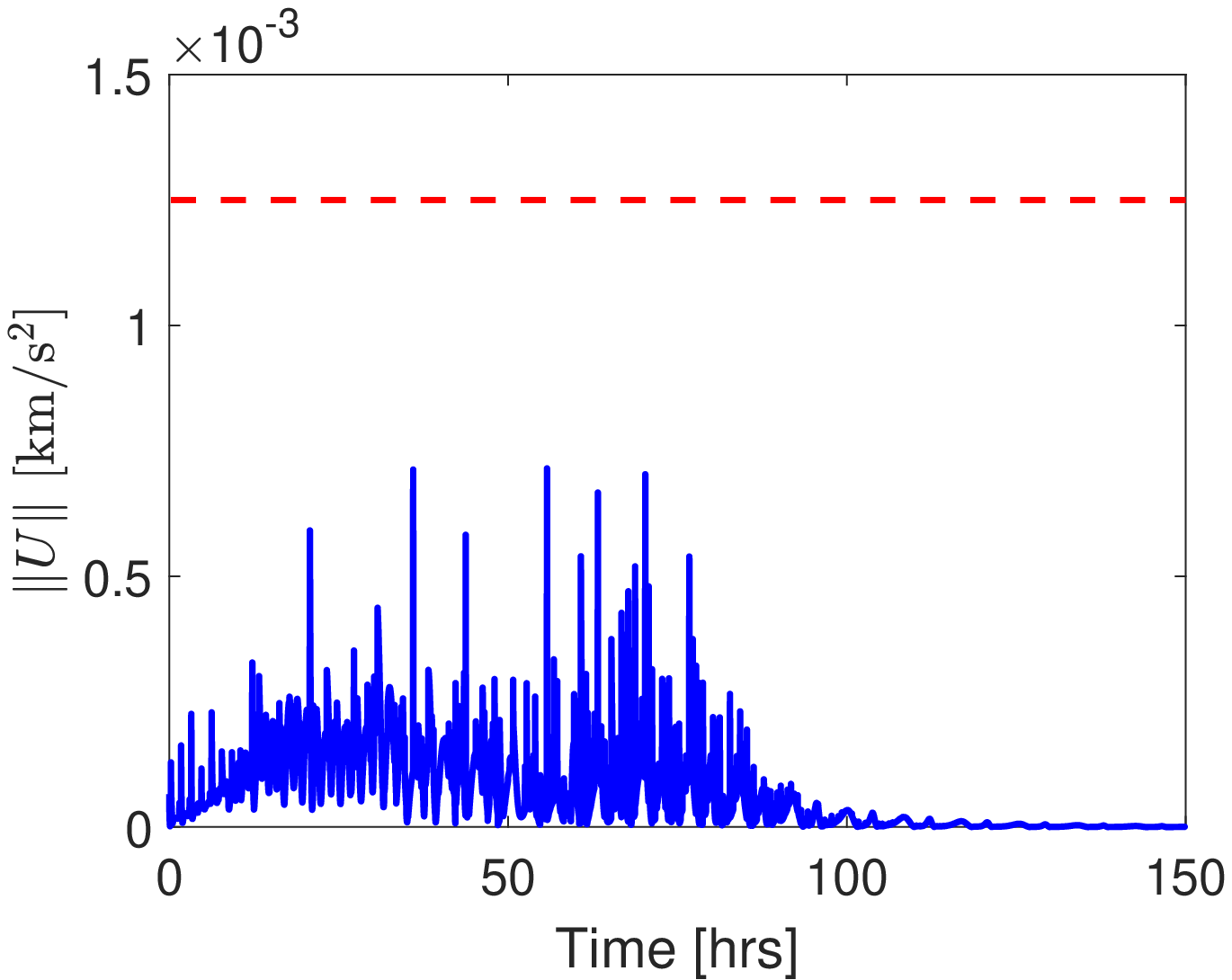}
}
\caption{Orbital transfer from a lower orbit to a higher orbit.  Left: The time histories of $r_p$, $r$ and $r_{\tt min}$ showing that the constraint (\ref{equ:maincon1}) is enforced. Right: The time histories of  $\|U\|$ (solid, blue) and of $U_{\tt max}$ (red, dashed)
showing that  the constraint (\ref{equ:maincon2}) is enforced.}\label{fig:Going up coss}
\end{figure}

\section{Varying Spacecraft Mass}\label{sec:7}
During the orbital transfer maneuvers, fuel is consumed by the thrusters causing the spacecraft mass to decrease over time according to
\begin{equation*}
    \dot m(t) = -\frac{F(t)}{I_{\tt sp} g_0},
\end{equation*}
where $m(t)$ is the mass of the spacecraft, $F(t)$ is the magnitude of the thrust force, $I_{\tt sp}$ [sec] is the thruster specific impulse and $g_0$ is the average Earth acceleration due to gravity at the sea level.

The constraint (\ref{equ:maincon2}) sets an upper limit on the spacecraft acceleration usable for the maneuver. The thrusters are limited by a maximum thrust force level, so, if the spacecraft mass is changing during the maneuver, the upper limit on the acceleration $U_{\tt max}$ in (\ref{equ:maincon2}) becomes time-varying.  In the implementation of m$^2$IRG, we set the upper acceleration limit at the time instant $t_k$ when $\tilde{X}(t_k)$ is decided on as
\begin{equation*}
   U_{\tt max}(t_k) = \frac{F_{\tt max}}{m(t_k)}.
\end{equation*}

In the simulations, we assumed that the spacecraft has Aerojet Rocketdyne MR-107T propulsion system\cite{rocketdyne}, with a maximum thrust of ($F_{\tt max}=0.125$ kN) the tank of which can take up to $39.39$ kg of hydrazine\cite{arianetank}. The spacecraft initial mass (dry mass plus fuel mass) was assumed to be $100$ kg, out of which $39.39$ kg was the fuel mass.

\begin{figure}[h]
\centering
\subfloat[]{
    \includegraphics[width=4.5cm]{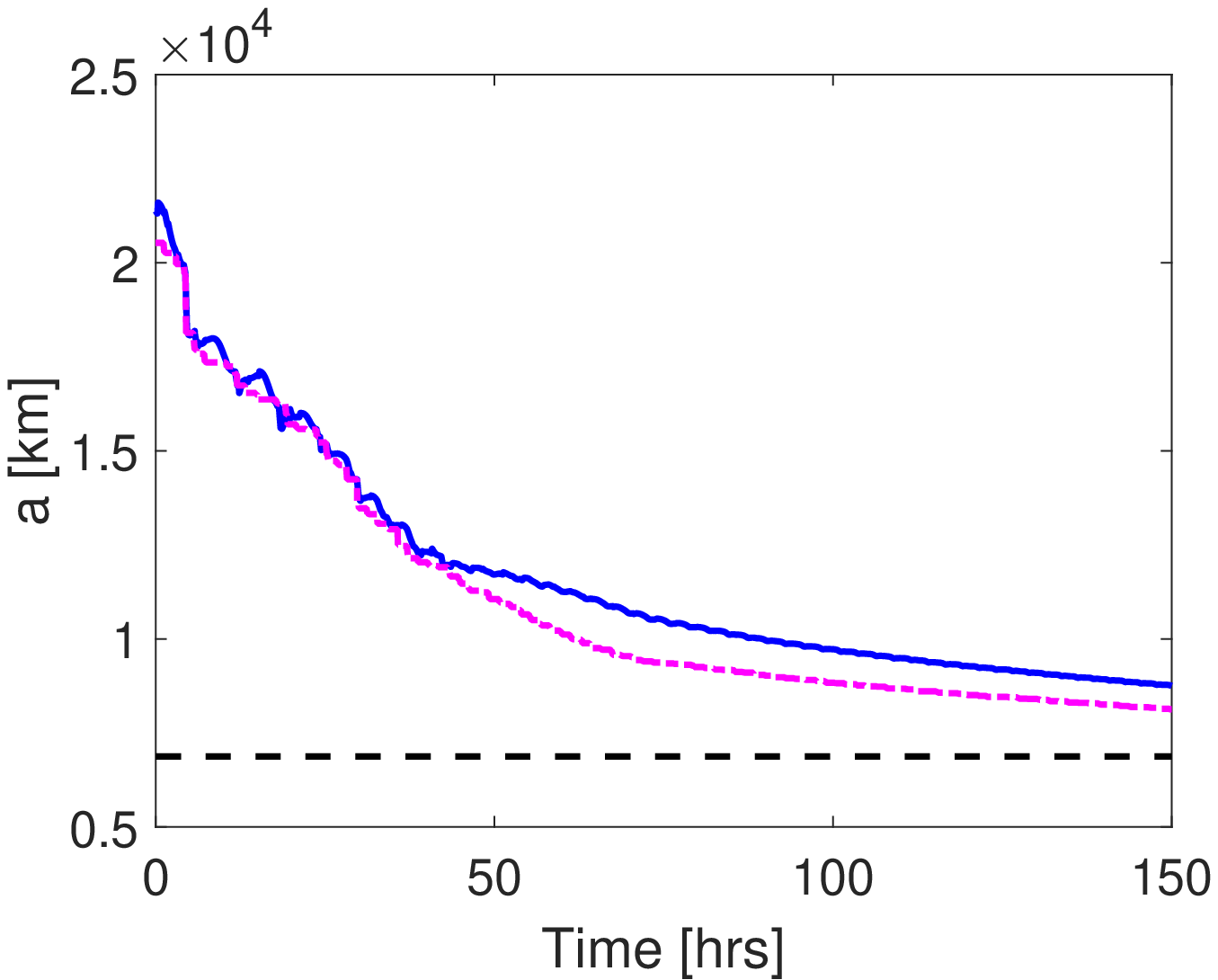}
}\qquad
\subfloat[]{
    \includegraphics[width=4.5cm]{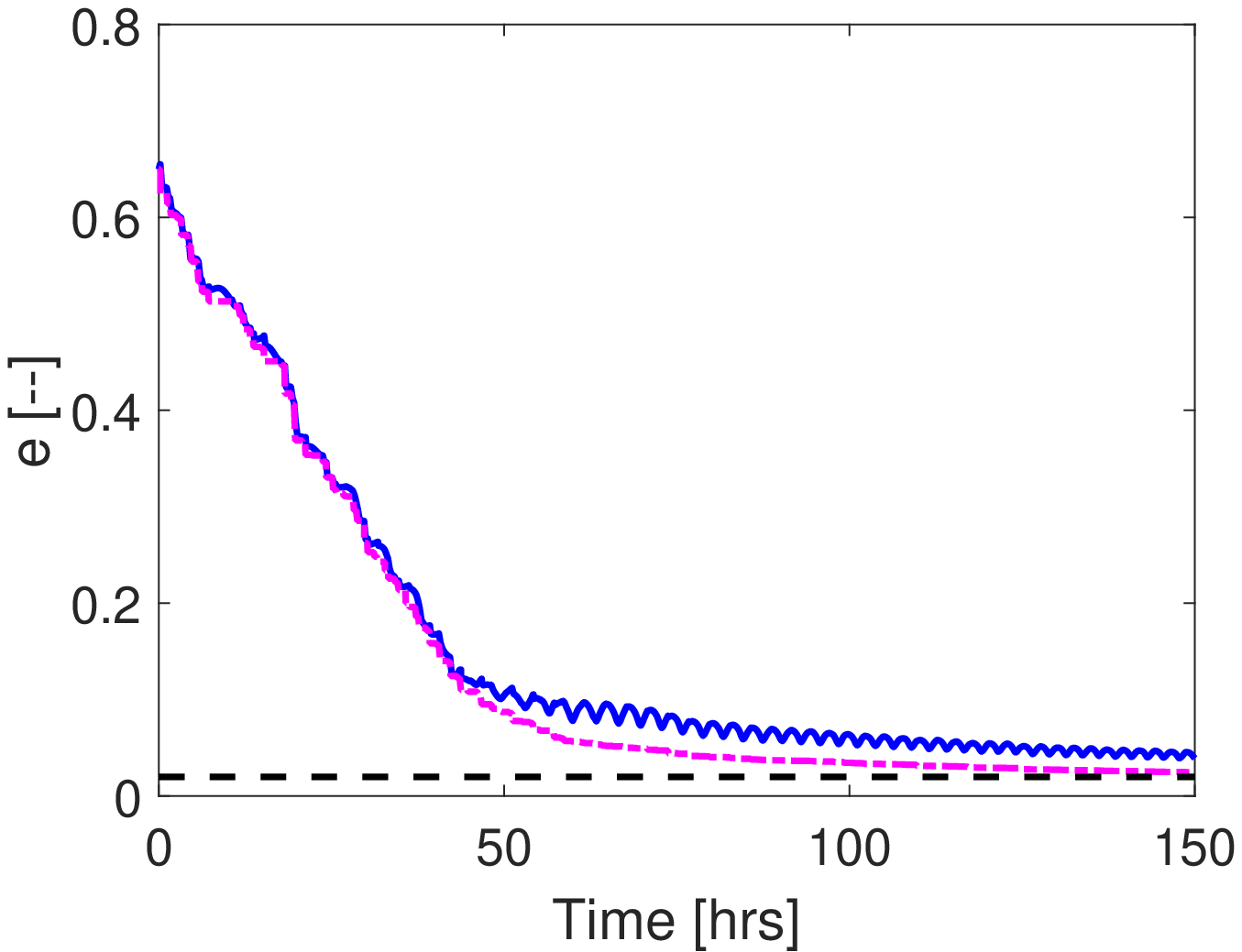}
}\qquad
\subfloat[]{
    \includegraphics[width=4.5cm]{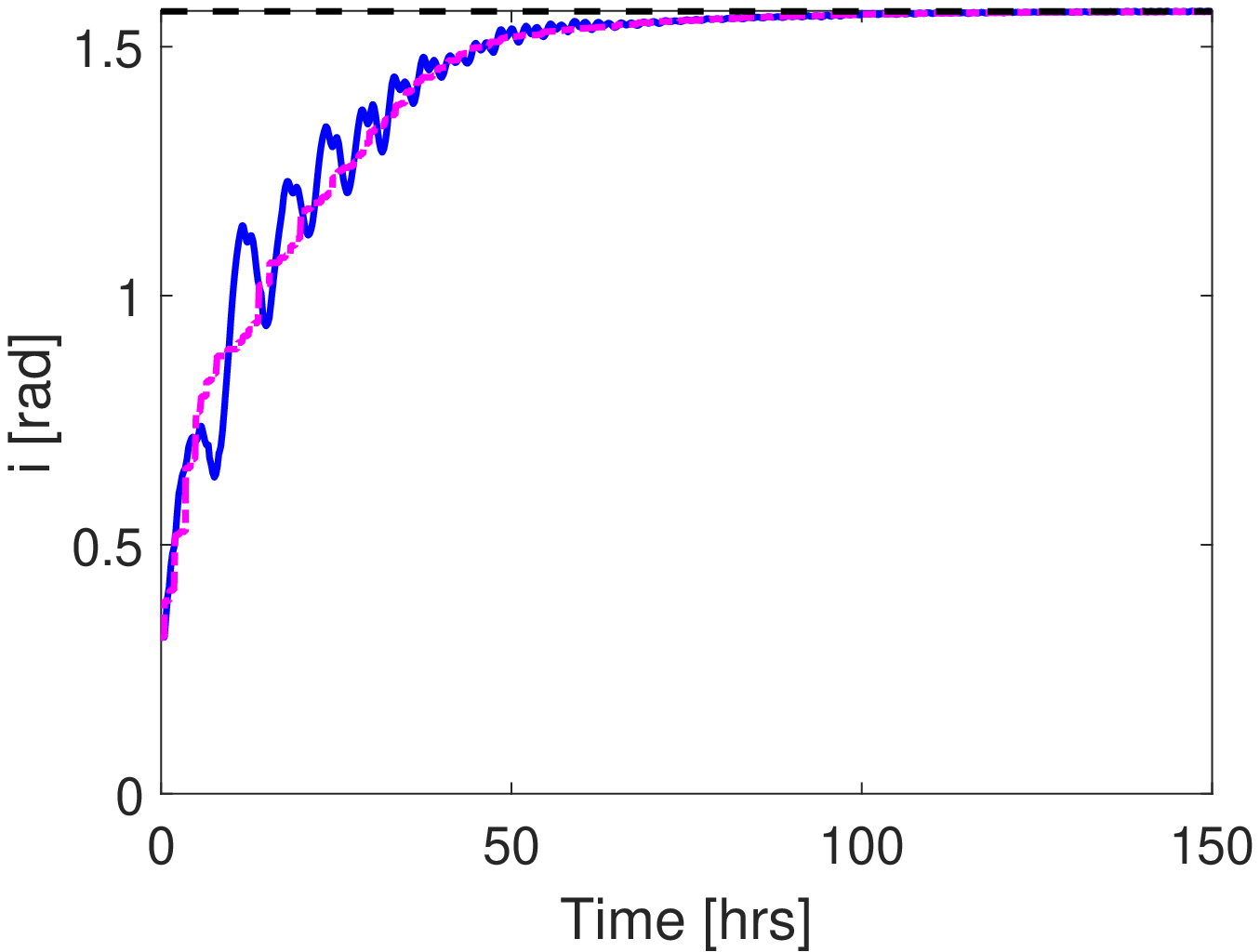}
}\vskip\baselineskip
\subfloat[]{
    \includegraphics[width=4.5cm]{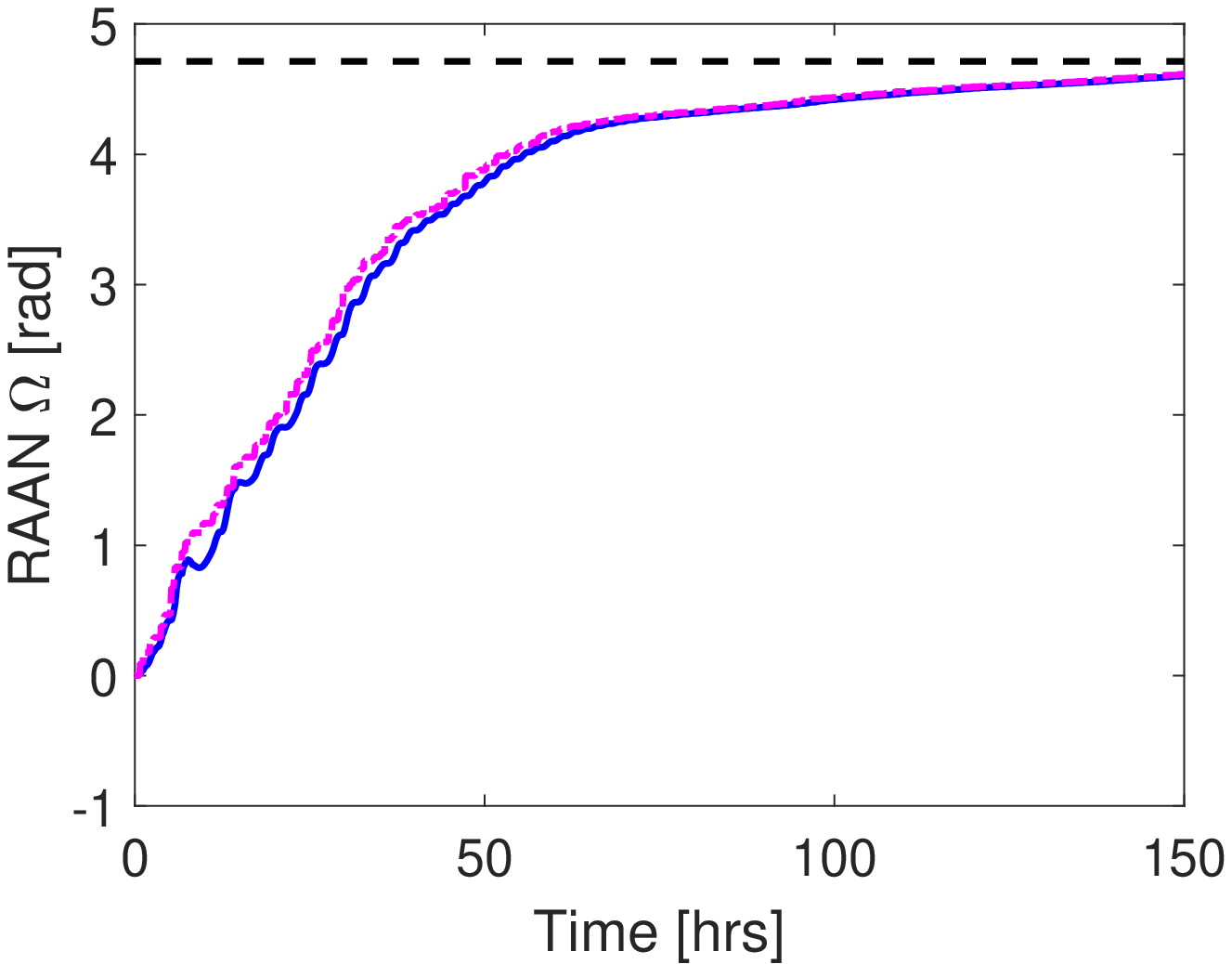}
}\qquad
\subfloat[]{
    \includegraphics[width=4.5cm]{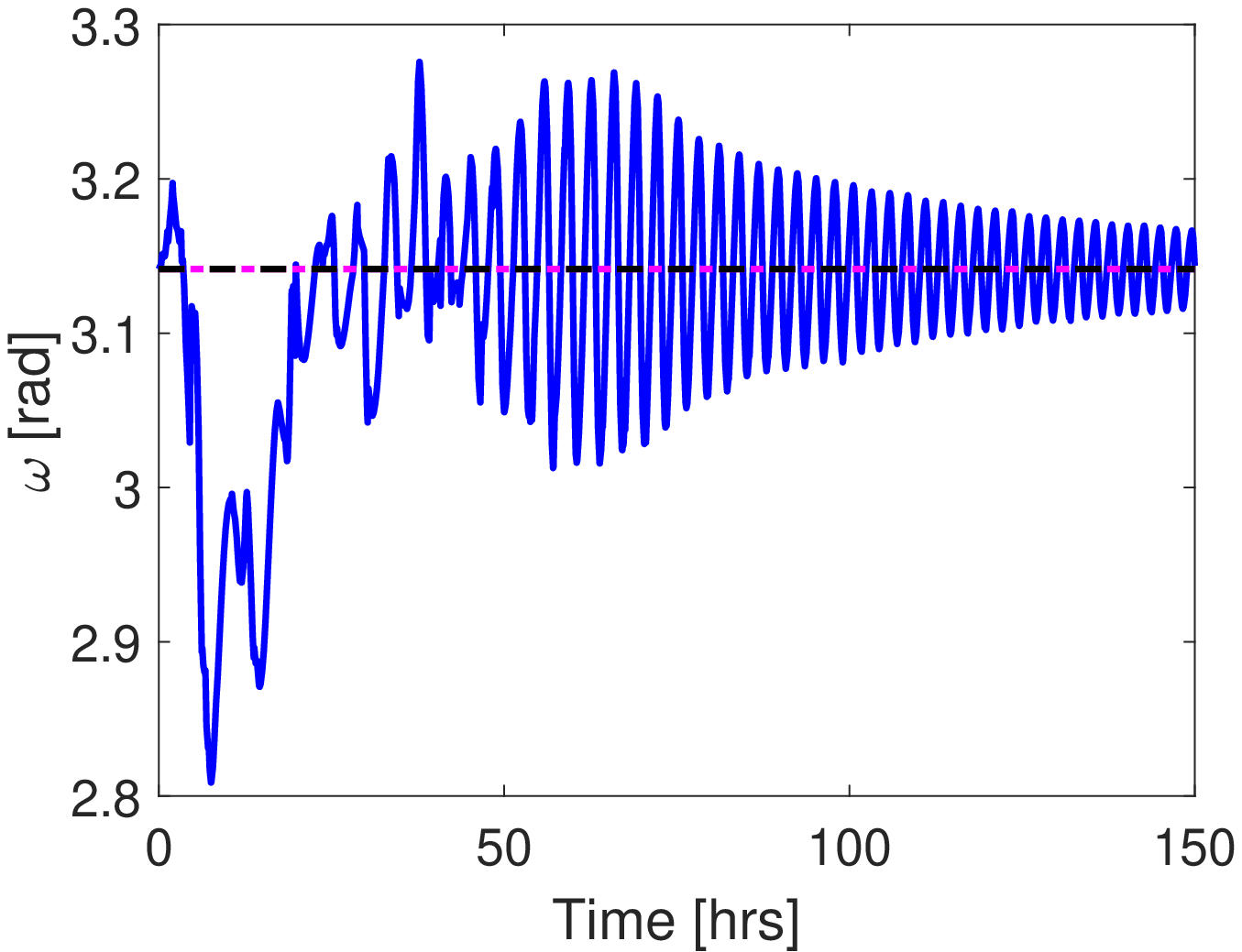}
}\qquad
\subfloat[]{
    \includegraphics[width=4.5cm]{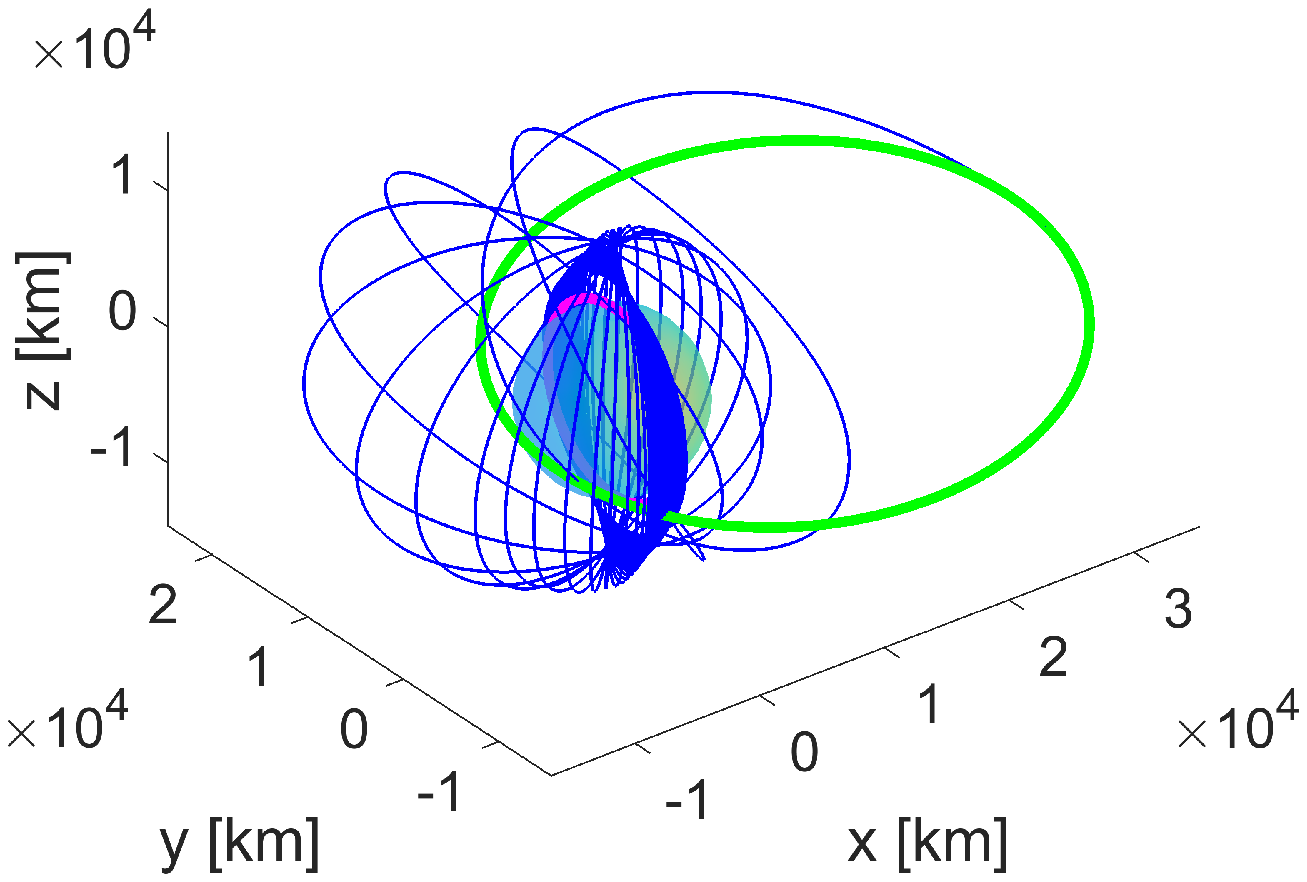}
    \label{subfig:Going down mass 3d}
}
\caption{Orbital transfer from higher orbit to lower orbit with m$^2$IRG when the spacecraft mass varies due to fuel consumption. The time histories of orbital elements evolution are shown in blue, compared  with the final desired state (black) and the reference generated by the m$^2$IRG (magenta). Figure~\ref{subfig:Going down mass 3d} shows the three dimensional trajectory followed by the spacecraft from the higher orbit (green) to the lower orbit (magenta).}\label{fig:going down mass}
\end{figure}

\begin{figure}[h]
\centering
\subfloat[]{
\includegraphics[width = 6 cm]{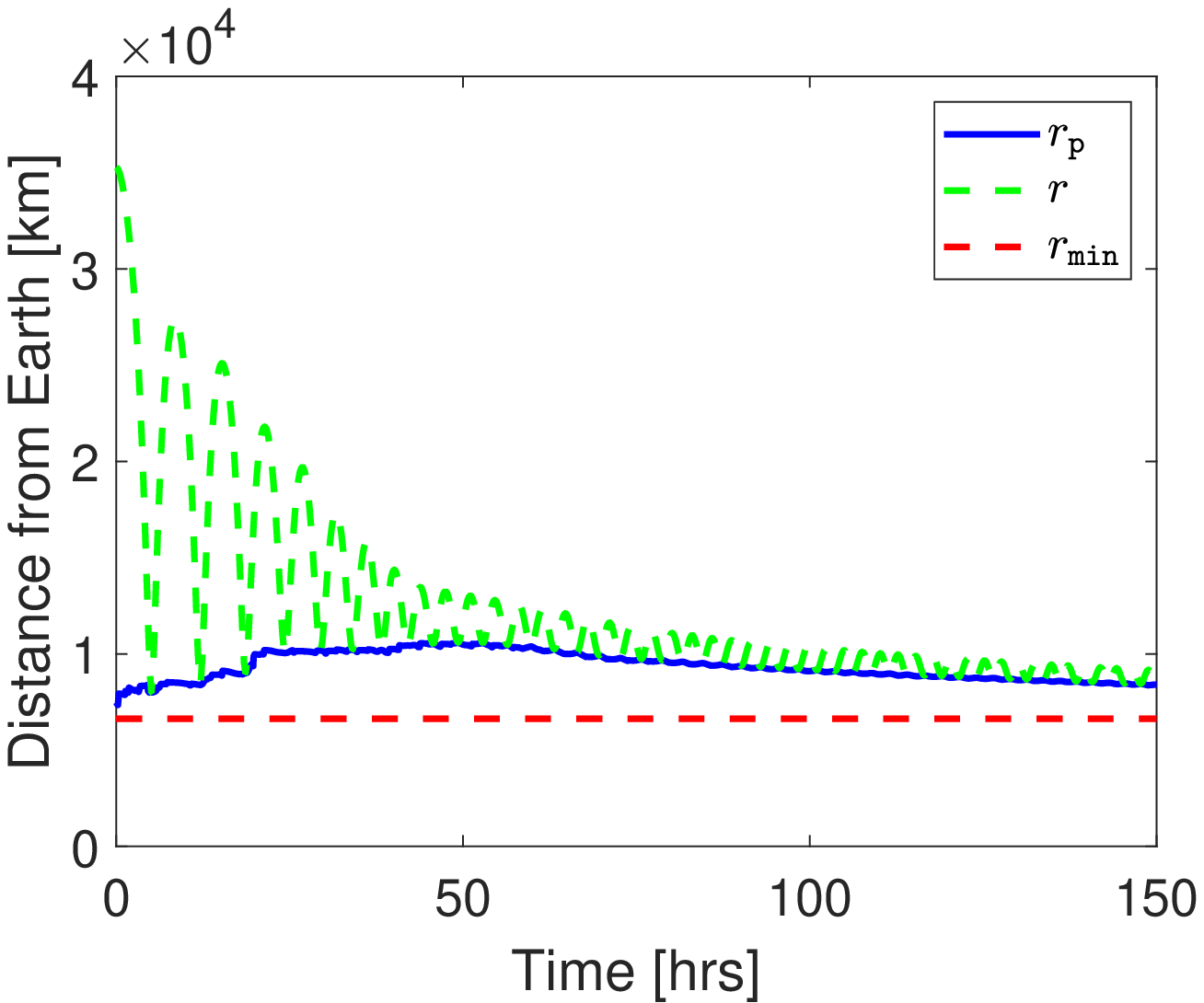}
}\qquad
\subfloat[]{
   \includegraphics[width=6cm]{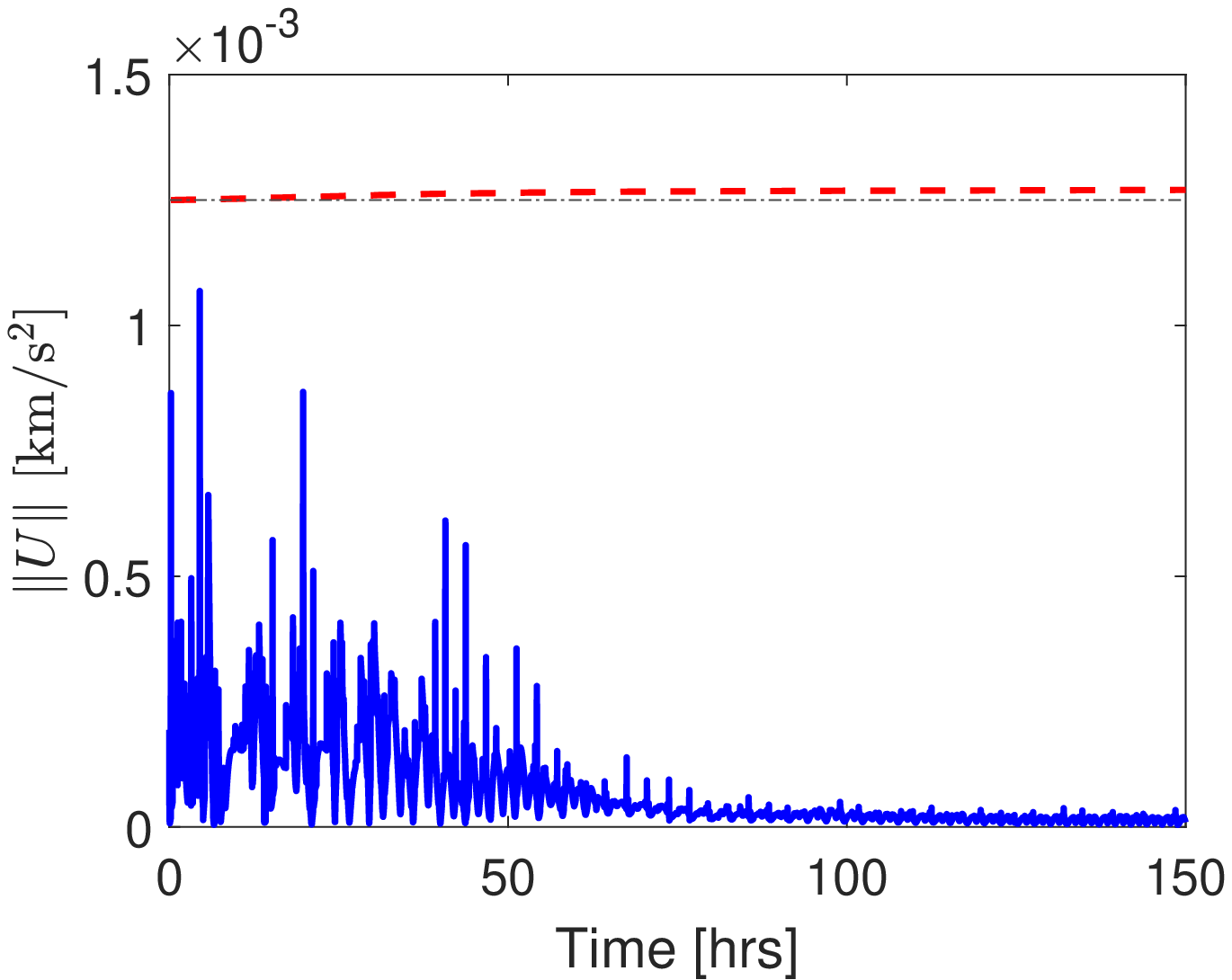}  
}
\caption{Orbital transfer from higher orbit to lower orbit with m$^2$IRG when the spacecraft mass varies due to fuel consumption.   Left: The time histories of $r_p$, $r$ and $r_{\tt min}$ showing that the constraint (\ref{equ:maincon1}) is enforced. Right: The time histories of  $\|U\|$ and of $U_{\tt max}$ showing the constraint (\ref{equ:maincon2}) is enforced.}\label{fig:mass variation}
\end{figure}

The time histories of the orbital elements for the transfer between a higher orbit to a lower orbit are shown in Figure~\ref{fig:going down mass}. Figure~\ref{fig:mass variation} is showing the more interesting result of changing the limitation on acceleration with time. The increase is small, due to reduced mass consumption rate - about 4\% difference from the initial satellite overall mass. The remarkable outcome, although, is that the scheme can handle also time-varying constraints.

\begin{figure}[h]
\centering
\includegraphics[width = 5 cm]{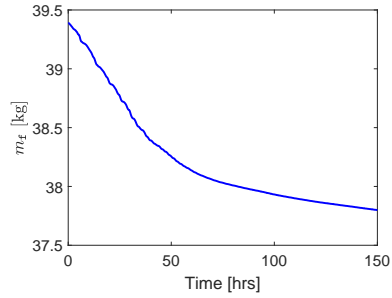}
\caption{Fuel mass change during the maneuver.}\label{fig:minrad mass var}
\end{figure}

\section{Online Prediction-Based Incremental Reference Governor}\label{sec:8}
A different approach to implementing m$^2$IRG involves changing the method for checking the acceptability of a given $\tilde{X}(t_k)$.  The use of sublevel sets of a Lyapunov function to bound the predicted trajectory and (\ref{equ:mycon1})-(\ref{equ:mycon3}) based on such sublevel sets to check for constraint violation can be 
replaced by an online prediction of the closed-loop spacecraft trajectory, through simulations over a sufficiently long horizon, and verifying if this predicted trajectory satisfies the constraints (\ref{equ:maincon1})-(\ref{equ:maincon3}).  This approach has been adopted in online prediction based reference governors\cite{bemporad1998reference,garone2017reference} and the initial variant of IRG\cite{tsourapas2009incremental}.  

\begin{figure}[h]
\centering
\subfloat[]{
    \includegraphics[width=4.5cm]{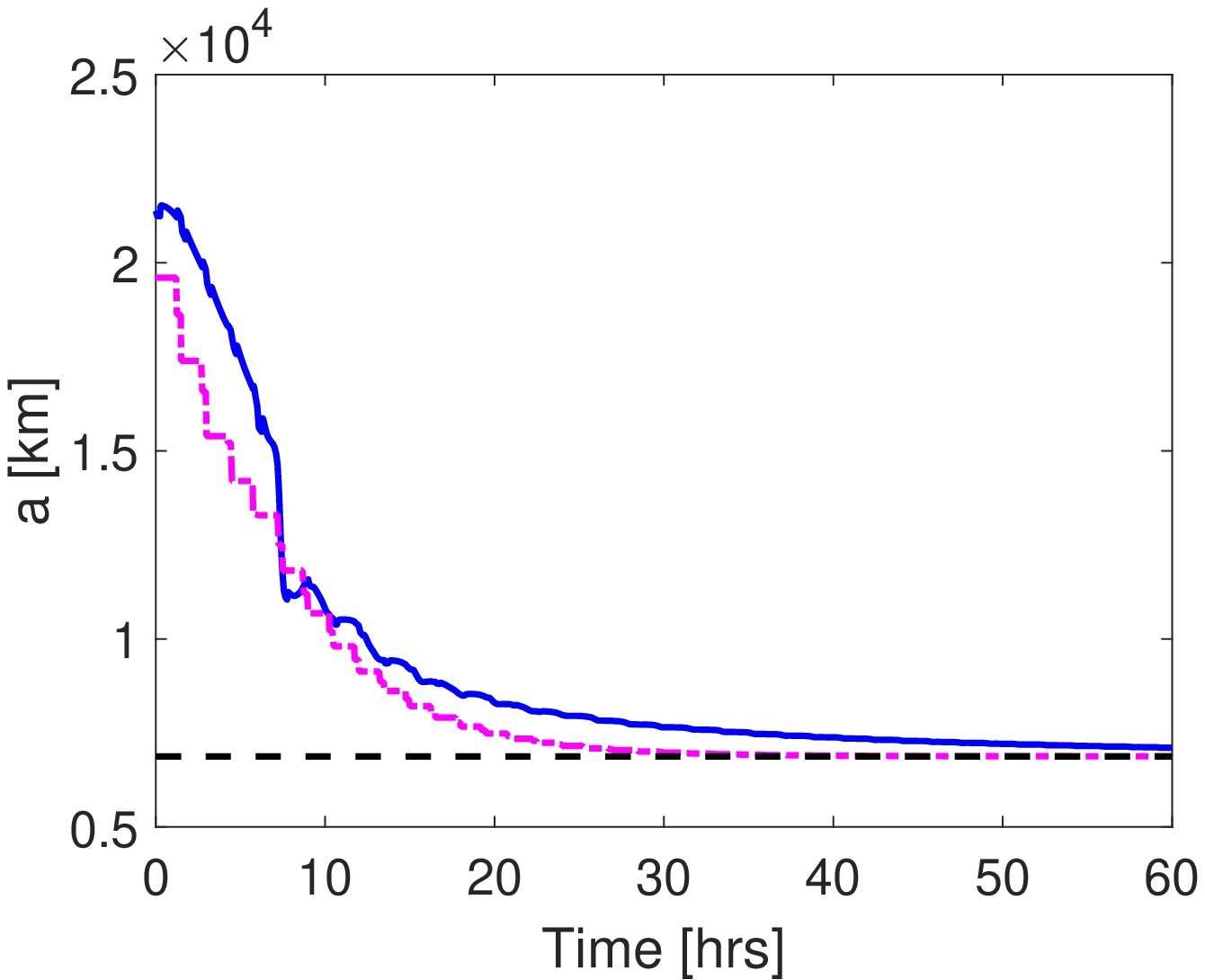}
}\qquad
\subfloat[]{
    \includegraphics[width=4.5cm]{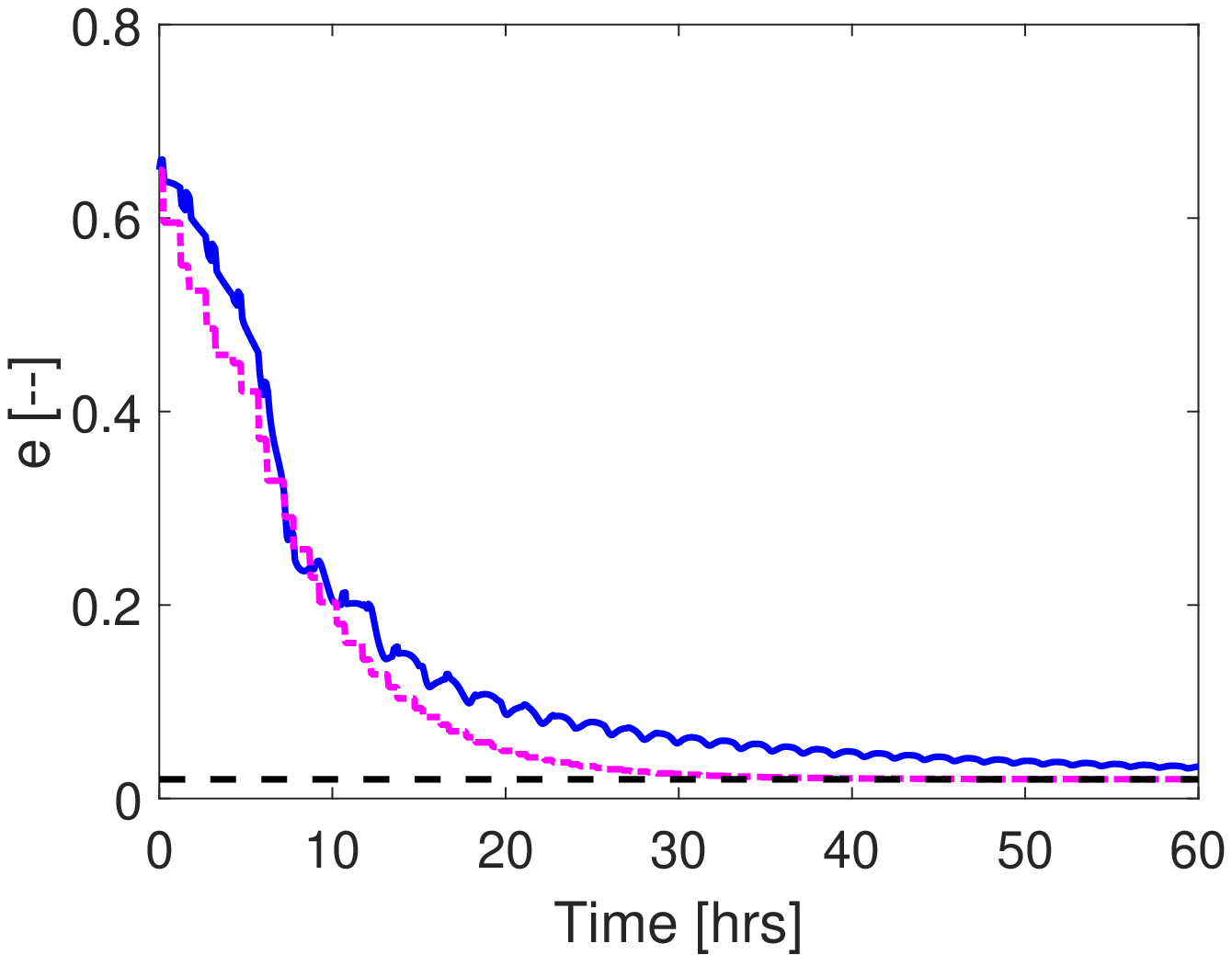}
}\qquad
\subfloat[]{
    \includegraphics[width=4.5cm]{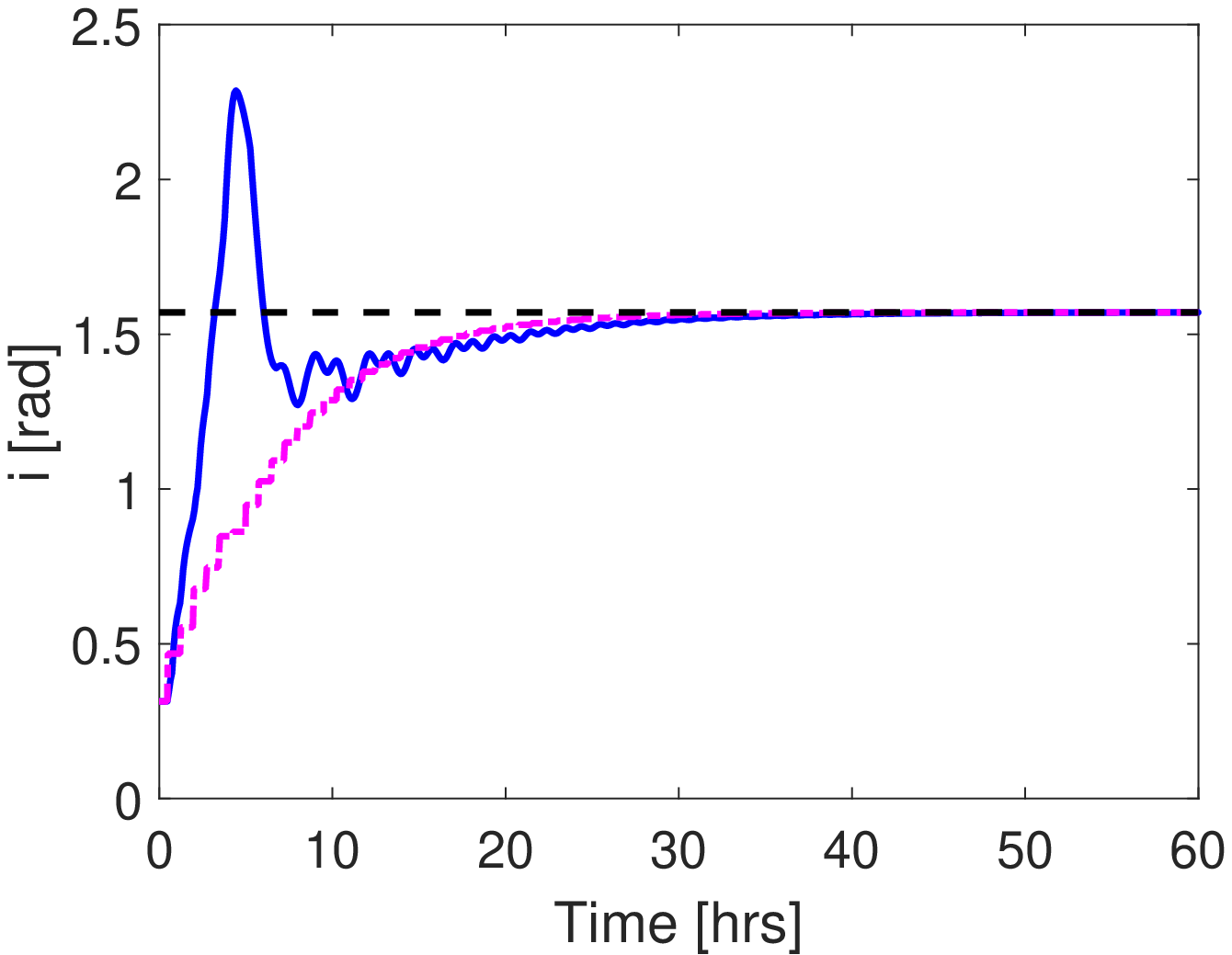}
}\vskip\baselineskip
\subfloat[]{
    \includegraphics[width=4.5cm]{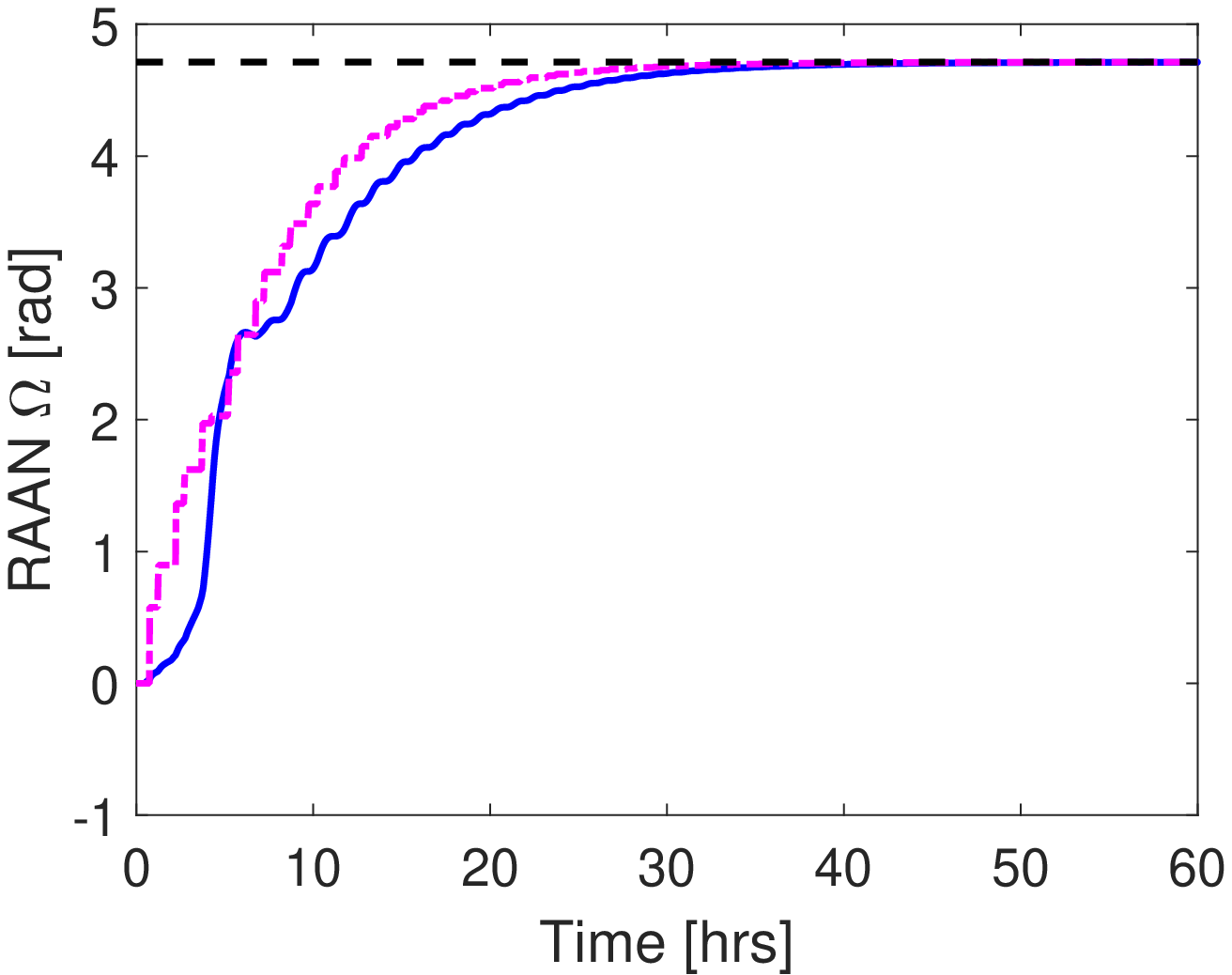}
}\qquad
\subfloat[]{
    \includegraphics[width=4.5cm]{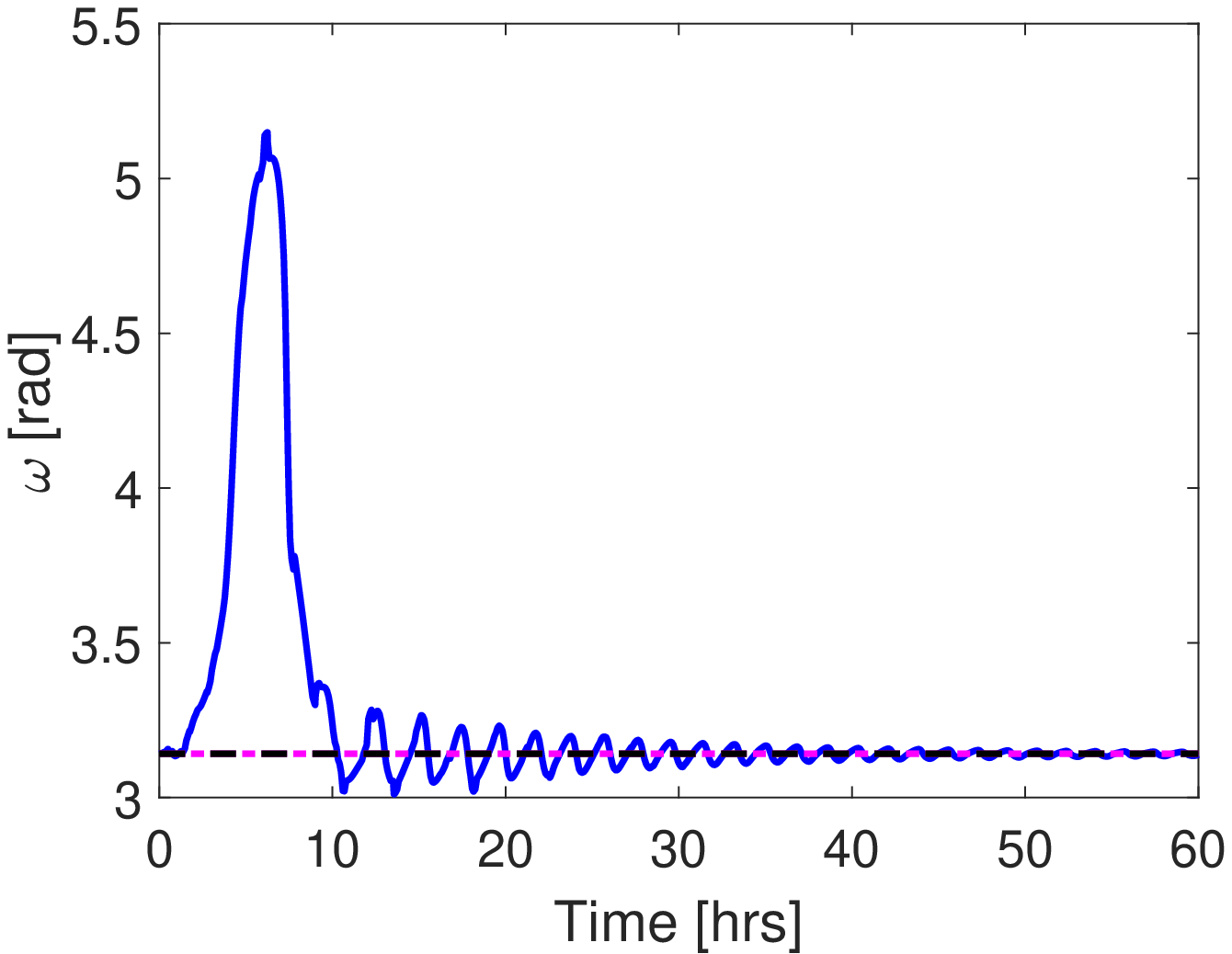}
}\qquad
\subfloat[]{
    \includegraphics[width=4.5cm]{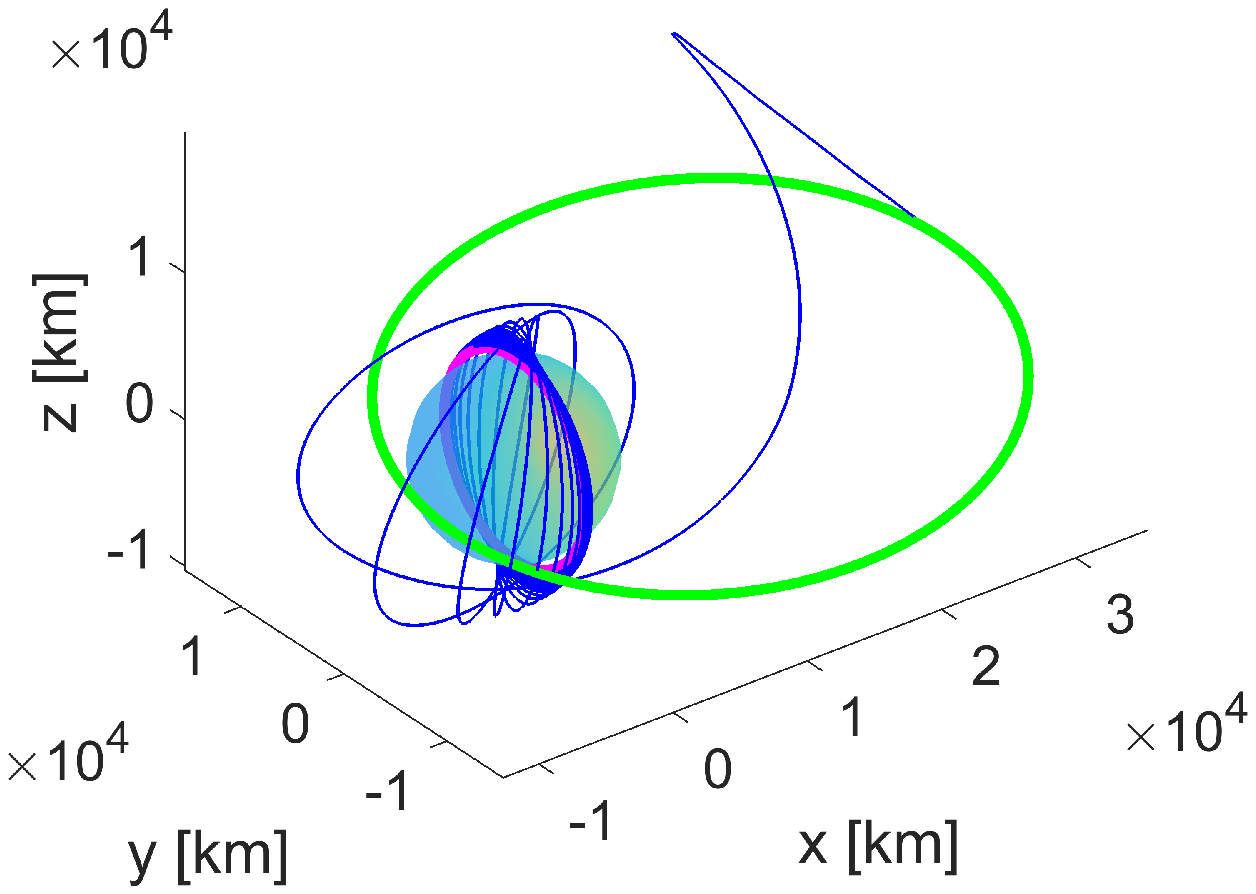}
    \label{subfig:PBIRG3d}
}
\caption{Orbital transfer from a higher orbit to a lower orbit with the online prediction-based m$^2$IRG { and with varying spacecraft mass}. The time histories of orbital elements  are shown in blue, compared  with the final desired state (dashed, black) and the reference generated by the m$^2$IRG (dash-dotted, magenta).  Figure~\ref{subfig:PBIRG3d} shows the three dimensional trajectory followed by the spacecraft from the initial (green) to the final (magenta) orbit.}\label{fig:PBIRGfigs}
\end{figure}

\begin{figure}[h]
\centering
\subfloat[]{
    \includegraphics[width=6cm]{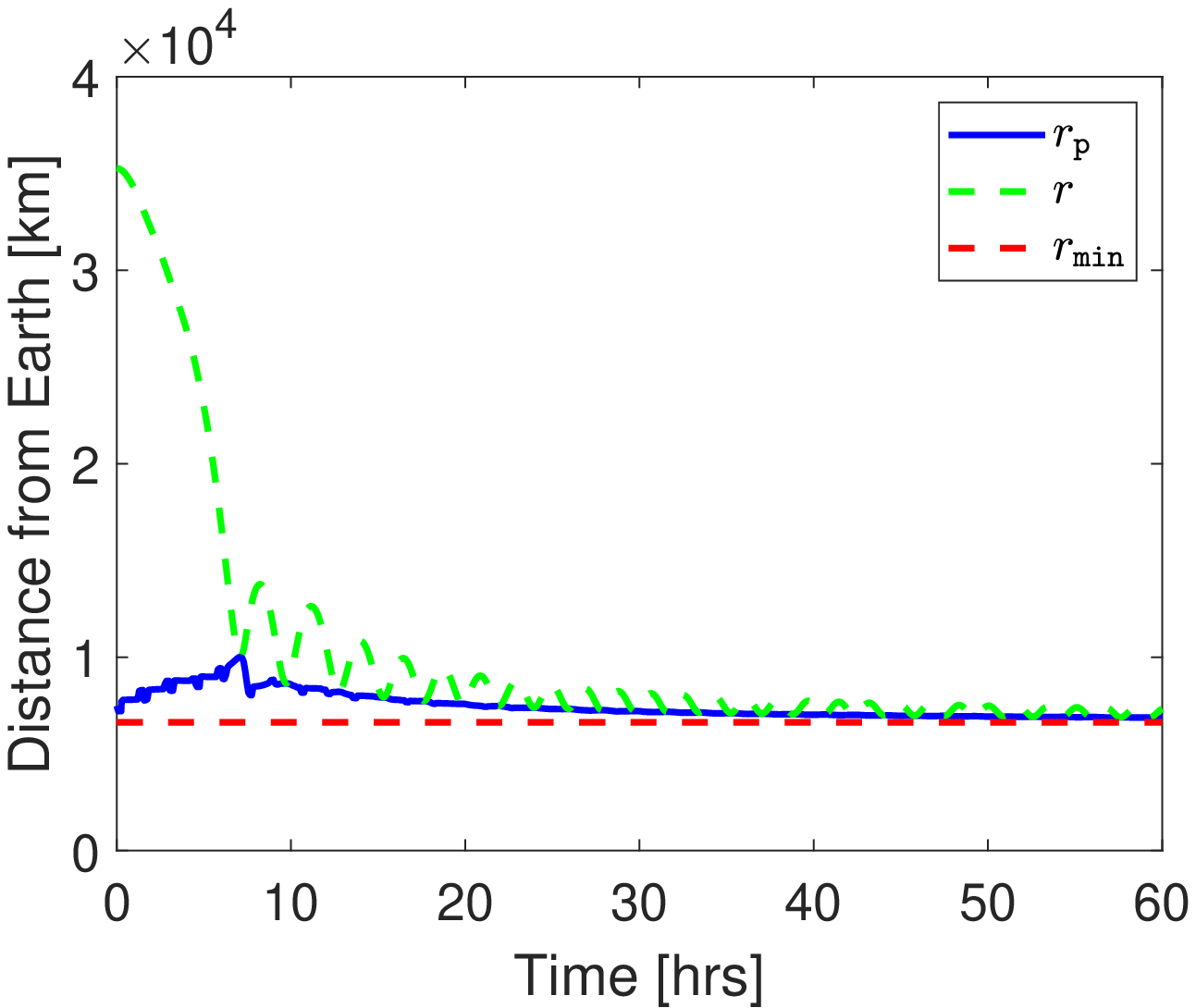}
}\qquad
\subfloat[]{
    \includegraphics[width=6cm]{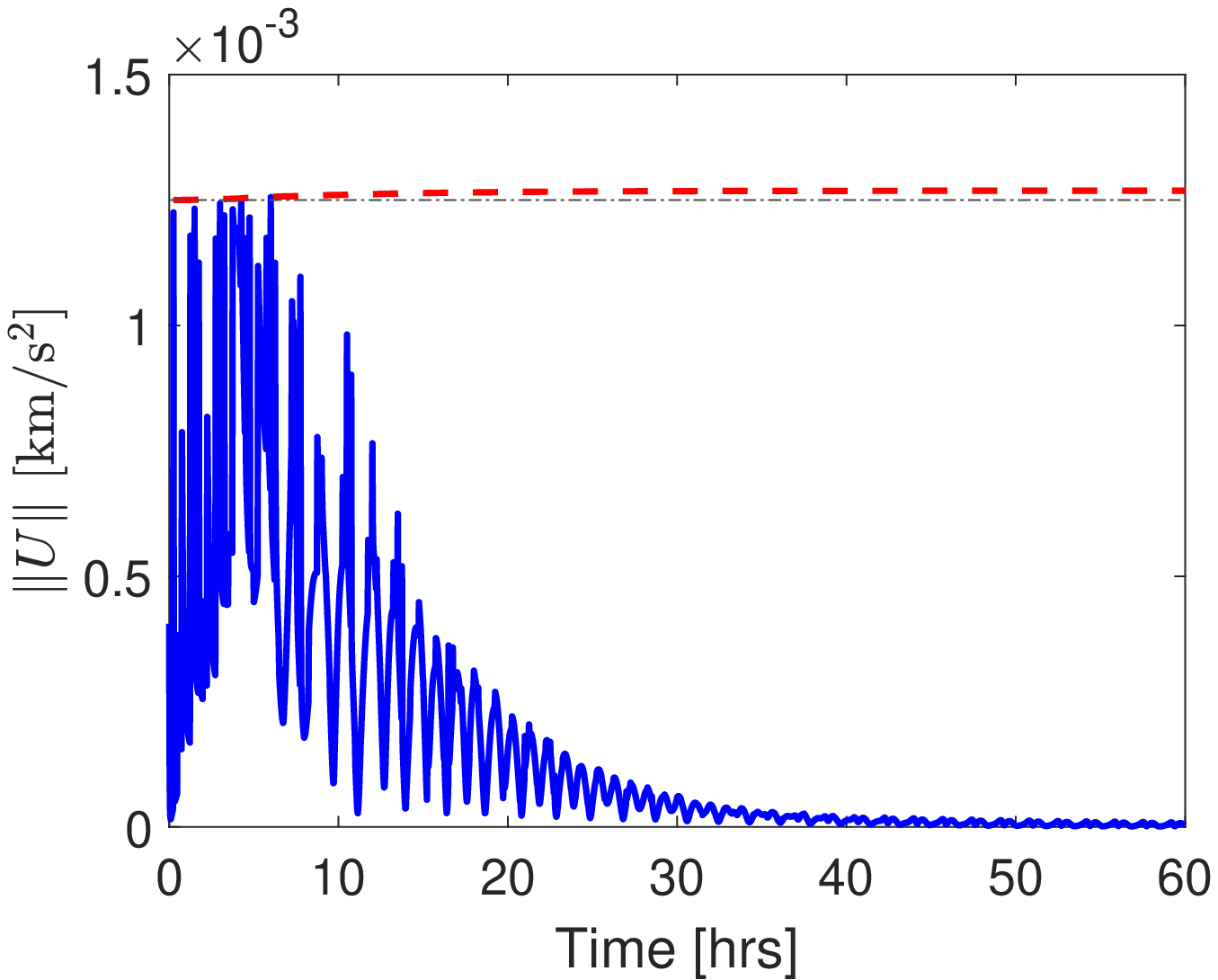}
}
\caption{ Orbital transfer from a higher orbit to a lower orbit with the online prediction-based m$^2$IRG { and with varying spacecraft mass}.  Left: The time histories of $r_p$, $r$ and $r_{\tt min}$ showing that the constraint (\ref{equ:maincon1}) is enforced. Right: The time histories of  $\|U\|$ (solid, blue) and of $U_{\tt max}$ (dashed, red)
showing  the constraint (\ref{equ:maincon2}) is enforced.}
\label{fig:PBIRG const}
\end{figure}

The simulations were performed with the same parameters as in Section~\ref{sec:7}. The prediction horizon was set to $10$ hours with constraints checked for violation every $20$ sec over the predicted trajectory. The time histories of the orbital elements and the three dimensional spacecraft trajectory are shown in Figure~\ref{fig:PBIRGfigs}. In combination with Figure~\ref{fig:PBIRG const}, it confirms that the constraints (\ref{equ:maincon1})-(\ref{equ:maincon3}) are satisfied.

\begin{figure}[h]
\centering
\includegraphics[width = 6.5 cm]{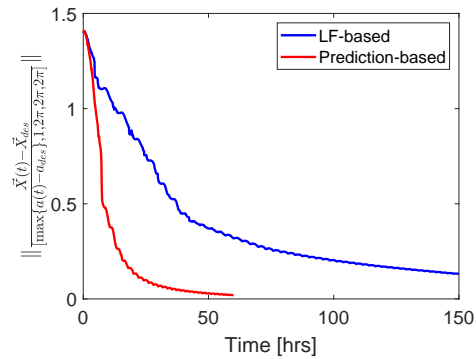}
\caption{Comparison between the convergence with Lyapunov Function (LF)-based m$^2$IRG and online prediction-based m$^2$IRG.}\label{fig:convcomp}
\end{figure}

As our simulation results show, the  online prediction based m$^2$IRG is less conservative and the closed-loop response in faster as compared to Lyapunov function based m$^2$IRG.  This is confirmed by Figure~\ref{fig:convcomp} which shows  the time history of the euclidean distance at time $t$ of the vector of orbital elements from the vector of their values on the target orbit, each normalized by their maximum value.  The normalization is performed to avoid bias induced by the significantly larger magnitude of the semi-major axis.

As compared to the Lyapunov function based m$^2$IRG, the computational footprint of the online prediction-based m$^2$IRG can actually be lower, as (\ref{equ:mycon1})-(\ref{equ:mycon3}) are solved online using numerical optimization while our implementation of the online prediction-based reference governor requires propagating forward EOMs, which is done using {\tt ode23s}. However, given that the convergence rate to the target equilibrium is determined by the ``richness'' of excitation provided by a particular trajectory, it is difficult to provide theoretical guarantees on the sufficient prediction horizon such that if the constraints are satisfied up to that horizon with constant $\tilde{X}$ and $P$, then they will be satisfied for all future times.  This issue is avoided by the Lyapunov function based m$^2$IRG that guarantees that the constraints are satisfied over the semi-infinite prediction horizon for the predicted trajectory if $\tilde{X}$ and $P$ are admissible.  

The knowledge of the Lyapunov function and the invariance properties of its sublevel sets can be exploited to implement the online prediction based m$^2$IRG with the variable prediction horizon (as in the conventional reference governor case~\cite{garone2017reference}).  In this approach, the constraints are checked up to a prediction horizon, $t_{\tt RG}$, which is the first time instant at which 
\begin{equation}\label{equ:predhor} 
V(X(t_k+t_{\tt RG}), \tilde{X}(t_k), P(t_k)) \leq \varepsilon,
\end{equation}
where $X(t_k+t_{\tt RG})$ denotes the predicted state and $\varepsilon>0$ is chosen sufficiently small so that all points in the set $\big\{X \in \mathbb{R}^5:~V(X,\tilde{X}(t_k),P(t_k)) \leq \varepsilon\big\}$ satisfy the constraints.   The estimates (\ref{equ:L1})-(\ref{equ:L2}) can in principle be exploited in estimating the required $t_{\tt RG}$ but their use is not straightforward as the contraction rate of Lyapunov function depends on the richness of excitation provided by a specific trajectory for which {\em a priori} analytical estimates are not available. Instead,  an additional criterion for deciding admissibility of  $\tilde{X}(t_k)$ and $P(t_k)$ in the m$^2$IRG algorithm can be added: If (\ref{equ:predhor}) is not satisfied for $0 \leq t_{\tt RG}
\leq t_{\tt RG}^{\tt max}$ for a sufficiently large $t_{\tt RG}^{\tt max}$, then
$\tilde{X}(t_k)$ and $P(t_k)$ are declared to be inadmissible in the m$^2$IRG algorithm at the time instant $t_k$.

\section{Concluding Remarks}\label{sec:9}
A multi-step multi-mode incremental reference governor (m$^2$IRG) has been developed to satisfy state and control constraints in feedback control of orbital transfer maneuvers.  Such m$^2$IRG augments a nominal Lyapunov controller which is derived based on Gauss Variational Equations of spacecraft dynamics.  The role played by persistence of excitation conditions has been highlighted in the nominal controller stability analysis.

Two variants of m$^2$IRG have been proposed.  The first variant relies on the invariance properties of the sublevel sets of the closed-loop Lyapunov function to ensure that constraints can be enforced over a semi-infinite prediction horizon.  The second variant relies on the online prediction of the closed-loop response over a finite and sufficiently long prediction horizon. Both m$^2$IRG variants adjust the reference command and accommodate switching of the Lyapunov controller gain.
Both m$^2$IRG variants have been able to successfully perform simulated orbital transfer maneuvers while satisfying the imposed state and control constraints and being robust to changes in the spacecraft mass due to fuel consumed during the maneuver.  The online prediction-based m$^2$IRG 
is less conservative in accommodating constraint violations as it relies on the actual trajectory predictions rather than their overbounds with the sublevel sets of the Lyapunov function. Consequently, the online predicting reference governor provides faster closed-loop response.

Analogs of GVEs have been derived\cite{battin1999introduction} for alternative sets of orbital elements, such as equinoctual orbital elements; these have a similar drift-free form as GVEs, however, are advantageous in term of not having a singularity at $e=0$ and $i=0$ and hence are suitable for handling circular and equatorial target orbits. Their treatment, while appears to be also tractable, is left as the subject for future research.

\section*{Acknowledgements}
The second author would like to acknowledge very useful discussions of the problem and insights into time-varying stabilization with Dr. Elena Panteley at L2S CNRS / Supelec, France.

\bibliography{ref.bib}
\end{document}